\documentclass{amsart}

\usepackage{amsmath,amssymb,amsthm}

\usepackage{pictexwd}

\DeclareMathOperator{\Ext}{Ext}
\DeclareMathOperator{\End}{End}

\DeclareMathOperator{\mo}{mod}
\DeclareMathOperator{\Hom}{Hom}
\DeclareMathOperator{\Ker}{Ker}
\DeclareMathOperator{\rad}{rad}
\DeclareMathOperator{\soc}{soc}
\DeclareMathOperator{\topp}{top}
\DeclareMathOperator{\bdim}{\bold {dim}}

\theoremstyle{plain}

\newcommand{\ssize}[1]{\smallmatrix #1\endsmallmatrix}
\newcommand{\arr}[2]{\arrow <1.5mm> [0.25,0.75] from #1 to #2}

\begin{document}

\title{Kronecker Modules Generated by Modules of Length 2}

\author{Claus Michael Ringel}

\email{ringel$@$math.uni-bielefeld.de}

\address{Fakult\"at f\"ur Mathematik, Universit\"at Bielefeld,
D-33501 Bielefeld, Germany\\ and
Department of Mathematics, Shanghai Jiao Tong University,
Shanghai 200240, PR China. \\}
\date{}
\maketitle

\begin{abstract}
Let $\Lambda$ be a ring and $\mathcal N$ a class of $\Lambda$-modules. A $\Lambda$-module is
said to be generated by $\mathcal N$ provided that it is a factor module of a direct sum
of modules in $\mathcal N$. The semi-simple $\Lambda$-modules are just the $\Lambda$-modules
which are generated by the $\Lambda$-modules of length 1. It seems that the modules
which are generated by the modules of length $2$ (we call them bristled modules)
have not attracted the interest they deserve.

In this paper we deal with the basic case of the Kronecker modules, these are
the (finite-dimensional) representations of an $n$-Kronecker quiver,
where $n$ is a natural number. We show that for $n\ge 3$,
there is an abundance of bristled
Kronecker modules.
\end{abstract}

\section{Introduction}

\subsection{}
Let $\Lambda$ be a ring and $\mathcal N$ a class of $\Lambda$-modules. A $\Lambda$-module is
said to be generated by $\mathcal N$
provided that it is a factor module of a direct sum
of modules in $\mathcal N$.
A $\Lambda$-module $M$ is said to be {\it $\mathcal N$-saturated}
provided that $\Ext^1(N,M) = 0$ for all modules $N\in \mathcal N$. Of course, {\it the
class of modules which are both generated by $\mathcal N$ and $\mathcal N$-saturated is
closed under extensions. In case $\Lambda$ is hereditary, the
class of modules which are both generated by $\mathcal N$ and $\mathcal N$-saturated is
also closed under factor modules.}
We consider the case where $\mathcal N$ is a class of
modules of length 2.
The modules which are generated by modules of length $2$ will be called
{\it bristled} modules; a {\it bristle} is by definition an
indecomposable module of length 2, and we denote by $\mathcal B$ the set of
isomorphism classes of bristles
(we should stress that a bristled module is not necessarily
the sum of its bristle submodules --- it is the sum of
submodules which are either bristles or simple).

Given any artin algebra $\Lambda$, we
denote by $\mo\Lambda$ the category of all (finitely generated) $\Lambda$-modules
and by $\tau$ the Auslander-Reiten translation in $\mo\Lambda$.

\subsection{}

We focus the attention to the $n$-Kronecker algebras and their modules,
the $n$-Kronecker modules (for convenience, we sometimes will refer to 
these algebras and their modules just as the Kronecker algebras and the Kronecker modules). 
    Despite of the importance of the Kronecker algebras,
    not much is known about Kronecker modules. Our approach to look at the bristled
    Kronecker modules may be seen as one of the first attempts
    to get insight into the structure of Kronecker modules in general.
    On the other hand, the Appendix B  will outline in which way the study
    of bristled Kronecker modules is an essential step for
    understanding bristled $\Lambda$-modules for any finite-dimensional $k$-algebra,
    where $k$ is an algebraically closed field.

Let $k$ be a field.
The $n$-Kronecker quiver $K(n)$ has two vertices, denoted by $1$ and $2$, and $n$ arrows
$\alpha_i\colon 1 \to 2$, with $1\le i \le n$.
$$
\hbox{\beginpicture
  \setcoordinatesystem units <3cm,3cm>
\put{$1$} at 0 0
\put{$2$} at 1 0
\arr{0.2 0.1}{0.8 0.1}
\arr{0.2 -.1}{0.8 -.1}
\put{$\alpha_1$} at 0.5 0.2
\put{$\alpha_n$} at 0.5 -.2
\put{$\vdots$} at  0.5 0.03
\endpicture}
$$
The $n$-Kronecker algebra $kK(n)$ is the path algebra of $K(n)$,
the $n$-Kronecker modules are the representations of $K(n)$ (over $k$).
Kronecker modules are written in the form $M =
(M_1,M_2;\alpha_1,\dots,\alpha_n)$,
where $M_1,M_2$ are two vector spaces and $\alpha_i\colon M_1\to M_2$ are
linear transformations, for $1\le i \le n$.

Note that $kK(n)$ is a hereditary artin algebra. 
Given any hereditary artin
algebra $\Lambda$, an indecomposable $\Lambda$-module is called {\it preinjective} provided
that it is of the form
$\tau^tI$ for some injective $\Lambda$-module $I$ and some $t\ge 0.$

An $n$-Kronecker module which is a bristle
is of the form $B(\lambda) = (k,k;\lambda_1,\dots,\lambda_n)$, where
$(\lambda_1,\dots,\lambda_n)$ is a non-zero element of $k^n$; the bristles
$B(\lambda)$ and $B(\lambda')$ are isomorphic if and only if
$(\lambda'_1,\dots,\lambda'_n)$ is a scalar multiple of $(\lambda_1,\dots,\lambda_n)$
(thus the set $\mathcal B$
of isomorphism classes of bristles may be identified with the projective
space $\mathbb P^{n-1}$).

\subsection{}\label{theorem}
      The essential part of the paper concerns the class of
      bristled Kronecker modules which are
      $\mathcal B$-saturated. As we have mentioned, this class of modules
      is closed under factor modules and extensions.
      As we will see, there is an abundance of such modules
      (but note that for $n\ge 2$,
        the bristles themselves do not belong to this class of modules, since
      they are not $\mathcal B$-saturated). Actually, it seems to be appropriate
      to select a suitable finite set $\mathcal B'$ of (isomorphism classes of) bristles 
      and to look at the (slightly smaller) class of modules
      which are generated by $\mathcal B'$ and which are $\mathcal B$-saturated (also this class is
      closed under factor modules and extensions). Here is the main theorem.
      	                    \medskip

{\bf Main Theorem.}
{\it Let $n \ge 3.$} {\it There is a set $\mathcal B_0$ of $n+2$
          bristles in $\mo kK(n)$ with the following properties:}
	        \smallskip
		            
{\rm(a)} {\it Any preinjective module is generated by
$\mathcal B_0$ and is $\mathcal B$-saturated.}
	        \smallskip
		            
{\rm(b)} {\it If $M$ is any module, then there is a number $t(M)$ such that all the modules
$\tau^tM$ with $t\ge t(M)$ are generated by
$\mathcal B_0$ and are $\mathcal B$-saturated.}
	        \medskip

For the choice of the set $\mathcal B_0$, see section \ref{choice}.
As we will see in the Addendum \ref{optimal}, the number $n+2$ is best possible:
{\it There is
no set of $n+1$ bristles which generates all the preinjective modules.}

As a consequence of Main Theorem, we see that for $n\ge 3$
there exists a regular module which generates all the preinjective modules, for
example the direct sum
of the modules in $\mathcal B_0$ (and there are also indecomposable regular modules
with this property, see \ref{ind}).       
In 1994, Kerner \cite{[K1]}
had asked whether such a module does exist.  In 2014, Kerner gave
a non-constructive proof for the existence of such a module \cite{[K2]}; here is now
an explicit example. 
Note that for $n = 2,$
any regular module generates only finitely many indecomposable preinjective
modules, see \ref{tame} and Appendix C.
      	                          \medskip

Using duality, we see that also the corresponding dual statements hold:
      \smallskip

(a)* {\it Any preprojective
module $M$ is cogenerated by
$\mathcal B_0$ and $\mathcal B$-cosaturated.}
	        \smallskip

(b)* {\it If $M$ is any module, then there is a number $t^*(M)$ such that all the modules
$\tau^{-t}M$ with $t\ge t^*(M)$ are cogenerated by
$\mathcal B_0$ and are $\mathcal B$-cosaturated.}
	        \medskip

\noindent
Of course, a module $M$ is said to be {\it cogenerated} by a class $\mathcal N$ of modules
provided that it is a submodule of a direct sum of modules in $\mathcal N$ (note that the modules
considered here are modules of finite length; otherwise we would have to take direct products
instead of direct sums). And $M$ is said to be
{\it $\mathcal N$-cosaturated} provided $\Ext^1(M,N) = 0$ for all
modules $N$ in $\mathcal N.$

\subsection{\bf  Outline of the paper.}

Section \ref{prelim} contains some preliminaries: we fix the
notation and recall some general properties of the categories $\mo kK(n)$.
Of course, the
structure of the categories $\mo kK(1)$ and $\mo kK(2)$ is well-known (whereas only
little is known about the categories $\mo kK(n)$ with $n \ge 3$).

The proof of theorem \ref{theorem} will be presented in
section \ref{proof} and some related results can be found in section \ref{related}.
In particular, we will show in section \ref{related}
that an indecomposable $\mathcal B$-saturated module
which is not simple is faithful. Section \ref{proof} (and most parts of
section \ref{related})
deal with $n$-Kronecker modules, where $n\ge 3$.
The cases $n = 1$ and $n = 2$ are discussed in section \ref{start}.

There are three appendices. The first is a short introduction to the paper \cite{[R1]}.
The set of bristle submodules of a Kronecker module $M$ without simple injective submodules
is in a natural way a projective variety. In \cite{[R1]} we show that any projective variety
occurs in this way. Appendix B deals with arbitrary artin algebras $\Lambda$
and their bristled modules. We are going to point out in which way the study of bristled
$\Lambda$-modules can be reduced to the study of bristled Kronecker modules.
As we have mentioned already, the existence of a regular module generating all
preinjective modules holds true for $K(n)$ with $n\ge 3$, but not for $K(2).$
It does not hold for any tame hereditary artin algebras, see the Appendix C.

\subsection{Acknowledgment.}
A weaker version of the Main Theorem
(as well as the result of \cite{[R1]})
were first presented at the conference in honor
of Jerzy Weyman's 60th birthday, April 2015, at the University of Connecticut.
The weaker version used a set of $2n\!-\!1$ bristles instead of the set $\mathcal B_0$
with the optimal cardinality $n+2$;
it was written up at a stay at the IPM in Isfahan, October 2016 and the author is grateful
to the Iranian hospitality.  The present version was finished in April 2017. 
The author is indebted to various mathematicians for questions and comments.

\section{Preliminaries}\label{prelim}

  The aim of this section is to recall
  some general properties of the categories
$\mo kK(n)$ of the $n$-Kronecker modules, see for example \cite{[R0]}.

Given any acyclic quiver $Q$ and a vertex $x$
of $Q$, we write $S(x)$, $P(x)$ and $I(x)$ for the corresponding indecomposable representation
of $Q$ which is simple, projective, or injective, respectively.
Thus, $S(x)_x = k$ and $S(x)_y = 0$ for any vertex $y\neq x$.  The
representation $P(x)$ is the projective cover, the representation $I(x)$ the injective
envelope of $S(x)$. 

\subsection{}There are two simple $n$-Kronecker modules,
namely the injective module $S(1) = I(1) = (k,0;0,\dots,0)$
and the projective module $S(2) = P(2) =
(0,k;0,\dots,0)$. 
										           
If $M = (M_1,M_2;\alpha_1,\dots,\alpha_n)$ is an $n$-Kronecker module,
we denote by $\bdim M = (\dim M_1,\dim M_2)$ (or also by $\binom{\dim M_1}{\dim M_2}$)
its {\it dimension vector}, it is an
element of the Grothendieck group $K_0(kK(n)) = \mathbb Z^2$.

On $K_0(kK(n))$, there is defined the bilinear form $\langle-,-\rangle =
\langle-,-\rangle_n$ via
$$
 \langle \binom ab,\binom{a'}{b'} \rangle = aa'+bb'-nab',
 $$
and one knows that
$$
 \langle\bdim M,\bdim M'\rangle = \dim \Hom(M,M') - \dim\Ext^1(M,M'),
 $$
for any pair $M,M'$ of $n$-Kronecker modules.

On $K_0(kK(n))$, the Coxeter transformation $\Phi\colon K_0(kK(n))
\to K_0(kK(n))$ is defined via $\Phi(a,b) = (n^2a-nb-a,na-b)$. It mimics the
Auslander-Reiten translation as follows: Let $M$ be an indecomposable
$n$-Kronecker module. If $M$ is non-projective, then
$$
 \bdim \tau M = \Phi(\bdim M).
 $$
 If $M$ is projective, then, of course, $\tau M = 0$
 (and we have $\Phi(\bdim P(i)) = - \bdim I(i)$
 for $i=1,2$).

For $n\ge 1$, a module is bristled if and only if it is
generated by the bristles and $S(2)$; and an indecomposable bristled module is either
isomorphic to $S(2)$ or else generated by the set $\mathcal B$ of bristles.

\subsection{The Auslander-Reiten quiver of $\mo kK(1)$.}
There are just three isomorphism classes of indecomposable 1-Kronecker modules,
namely $S(2) = P(2), B = P(1) = I(2)$ and $I(1) = S(1)$, and the Auslander-Reiten quiver
of $\mo kK(1)$ is of the form
$$
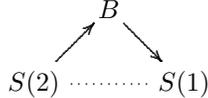
\hbox{\beginpicture
  \setcoordinatesystem units <1cm,1cm>
\put{$S(2)$} at 0 0
\put{$B$} at 1 1
\put{$S(1)$} at 2 0
\arr{0.3 0.3}{0.8 0.8}
\arr{1.2 0.8}{1.7 0.3}
\setdots <1mm>
\plot 0.5 0  1.5 0 /
\endpicture}
$$
(we write $B$ instead of $P(1) = I(2)$ in order to stress that this module
is a bristle; of course, it is the only bristle).

\subsection{The Auslander-Reiten quiver of $\mo kK(n)$ with
$n \ge 2$.} We recall that the indecomposable representations of
$K(n)$ fall into three classes: there are the preprojective ones, the regular ones
and the preinjective ones.

There are countably many (isomorphism classes of) indecomposable preinjective modules,
we denote these modules by $I_t$ with $t\in \mathbb N_0$, as follows:
$I_0 = S(1) = I(1)$,
$I_1 = I(2)$,  and $I_t = \tau I_{t-2}$
for $t\ge 2$, thus they are ordered in such a way that
$\Hom(I_t,I_{t'}) = 0$ if and only if $t < t'$.

Similarly, there are the indecomposable preprojective modules $P_0 = S(2) = P(2), P_1 = P(1),
P_2, P_3, \dots$, but this notation will not be needed in an essential way in the paper.

The Auslander-Reiten quiver looks as follows:
$$
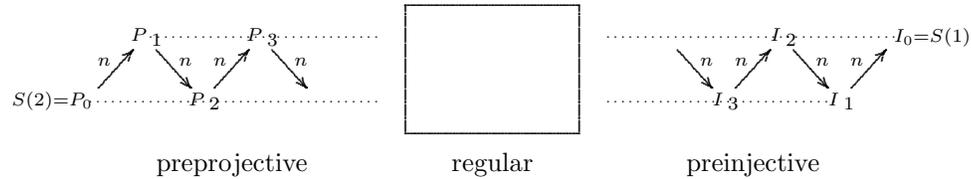
\hbox{\beginpicture
  \setcoordinatesystem units <.77cm,.85cm>
  \put{preprojective} at -4.5 -1
\put{regular} at 0 -1
\put{preinjective} at 4.5 -1
\put{$\ssize {S(2) =P_0}$} at -7.6 0
\put{$\ssize {I_0= S(1)}$} at 7.6 1
\put{$\ssize I_1$} at 6 0
\put{$\ssize I_2$} at 5 1
\put{$\ssize I_3$} at 4 0
\put{$\ssize P_1$} at -6 1
\put{$\ssize P_2$} at -5 0
\put{$\ssize P_3$} at -4 1
\plot -1.5 -.5  1.5 -.5  1.5 1.5  -1.5 1.5  -1.5 -0.5 /
\arr{3.2 0.8}{3.8 0.2}
\arr{4.2 0.2}{4.8 0.8}
\arr{5.2 0.8}{5.8 0.2}
\arr{6.2 0.2}{6.8 0.8}
\arr{-3.8 0.8}{-3.2 0.2}
\arr{-4.8 0.2}{-4.2 0.8}
\arr{-5.8 0.8}{-5.2 0.2}
\arr{-6.8 0.2}{-6.2 0.8}

\setdots <1mm>
\plot 2 1  7 1 /
\plot 2 0  6 0 /
\plot -2 1  -6 1 /
\plot -2 0  -7 0 /
\multiput{} at 0 1.8  0 -1.4 /
\multiput{$\ssize n$} at -6.7 0.6  -5.3 0.6  -4.7 0.6  -3.3 0.6
  6.3 0.6  5.7 0.6  4.3 0.6  3.7 0.6 /
  \endpicture}
  $$
  where we have endowed the arrows by the number $n$ in order to indicate that
  the corresponding bimodule of irreducible maps is $n$-dimensional.
  Note that there are no non-zero maps from a regular module to a
  preprojective module, or from a preinjective module to a
  preprojective or a regular module. All non-zero
  maps in the preprojective component and in the preinjective component go from
  left to right. The components in the regular part are of type $\mathbb A_\infty$.
  For $n = 2$, they are regular tubes; for $n\ge 3$ they are of the form $\mathbb Z A_\infty$.

It is important to be aware that the bristles are regular. Thus, any
indecomposable bristled module different from $S(2)$ is generated by bristles and
therefore is regular or preinjective.

\section{The bristled modules for $n=1$ and $n = 2.$}\label{start}

\subsection{The case $n = 1$.}
As we have mentioned, there is (up to isomorphism) a unique bristle, namely
$B = P(1) = I(2)$;
this bristle is projective and injective.
All the $1$-Kronecker modules are
$\mathcal B$-saturated bristled modules.

\subsection{The case $n = 2$.}

Recall that we have indexed the (isomorphism classes of the) bristles
by the elements of the projective space $\mathbb P^{n-1} = \mathbb P^1.$
For $c \in k$, let  $B_c = B(1:c)$,
and let $B_\infty = B(0:1)$.
    \medskip

{\bf Proposition.} {\it Let $n = 2$. Let $J$ be a subset of $k\cup\{\infty\}$.
Then $I_t$ is generated by the bristles $B_c$ with $c \in J$
if and only if $|J|  \ge  t+1.$}

\begin{proof}
Note that the dimension
vector of $I_t$ is $(t+1,t)$.
First, assume that
$I_t$ is generated by the bristles $B_c$ with $c \in J$.
Since $\dim\Hom(B_c,M) = 1$ for all bristles $B_c$,
it follows that $I_t$ is a factor module of the module
$\bigoplus_{c \in J} B_c$ with dimension vector $(|J|,|J|),$
therefore $|J| \ge t+1.$

On the other hand,
let $I_t = (M_1,M_2;\alpha_1,\alpha_2)$.
We can assume that $M_1$ has the basis $e_0,e_1,\dots, e_t$, that
$M_2$ has the basis $e'_1,\dots,e'_t$, that $\alpha_1(e_i) = e'_{i+1}$ for 
$0 \le i < t$ and $\alpha_1(e_t) = 0$ and that 
$\alpha_2(e_0) = 0$ and $\alpha_2(e_i) = e'_i$
for $1\le i \le n$, thus $M$ may be visualized as follows:
$$
\hbox{\beginpicture
  \setcoordinatesystem units <.8cm,1cm>
\put{$e_0$} at 0 1 
\put{$e_1$} at 2 1 
\put{$e_2$} at 4 1 
\put{$e_{t-1}$} at 7 1 
\put{$e_t$} at 9 1 
\put{$e'_1$} at 1 0 
\put{$e'_2$} at 3 0 
\put{$e'_t$} at 8 0 
\put{$e_0$} at 0 1 
\put{$e_0$} at 0 1 
\arr{0.3 0.7}{0.7 0.3}
\arr{1.7 0.7}{1.3 0.3}
\arr{2.3 0.7}{2.7 0.3}
\arr{3.7 0.7}{3.3 0.3}
\arr{7.3 0.7}{7.7 0.3}
\arr{8.7 0.7}{8.3 0.3}
\multiput{$\cdots$} at 5.5 0.9  5.5 0.45  5.5 0 /
\multiput{$\ssize\alpha_1$} at 0.7 0.7  2.7 .7  7.7 .7 /
\multiput{$\ssize\alpha_2$} at  1.3 0.7  3.3 .7  8.3 .7 /
\endpicture}
$$
For any $c\in k$, the element $m_c = \sum_{i=0}^t c^i e_i$ in $M_1$
is sent under $\alpha_1$ to $\alpha_1(m_c) = 
\sum_{i=1}^t c^{i-1}e'_i$ and under $\alpha_2$
to  $\alpha_2(m_c) = \sum_{i=1}^t c^{i}e'_i 
= c\cdot\alpha_1(m_c). $
This shows that $m_c$ generates a bristle of the form $B_c.$
The Vandermonde determinant shows that for $t+1$ pairwise different elements 
$c$
in $k$, the corresponding elements $m_c$ in $M_1$
are linearly independent, thus form a basis. But this implies that $M$ is a
bristled module.

There is one additional submodule of $M$ which is a bristle, namely the submodule
generated by $m_\infty = e_t$; it is isomorphic to $B_\infty.$ Given
$t$ pairwise different scalars 
$c_1,\dots,c_t$ in $k,$ then also 
the elements $m_{c_1},\dots,m_{c_t}, m_\infty$
are linearly independent: the determinant of the
$(t+1)\times(t+1)$ matrix whose rows are the vectors
$m_{c_1},\dots,m_{c_t}, m_\infty$ is the product of the
$t\times t$ Vandermonde determinant for the scalars   
$c_1,\dots,c_t$ and the determinant of the $1\times 1$-matrix with entry $1$,
thus non-zero.
\end{proof}

\subsection{} {\bf Corollary.}
{\it Let $n =2$ and $q = |k|.$
An indecomposable module $M$ different from $S(2)$ is bristled
if and only if $M$ is one of the bristles $B_c$ with $c\in k\cup\{\infty\}$
or $M = I_t$ is preinjective with $0\le t \le q$.

The preinjective bristled  
modules are $\mathcal B$-saturated,
the module $S(2)$ and the bristles are not $\mathcal B$-saturated.}

\begin{proof}
Let $M$ be an indecomposable bristled module.
First, assume that $M$ is regular. If the regular socle $U$ of $M$
is a bristle $B$, then $\dim\Hom(B,M) = 1$ and $\Hom(B',M) = 0$
for any bristle not isomorphic to $B$, thus $M = U = B$. Or else $\Hom(B,M) = 0$
for any bristle $B$, thus $M = 0$ (since $M$ is generated by bristles),
a contradiction.
This shows that the only indecomposable regular bristled modules are the bristles.

It remains to decide whether the bristled modules are $\mathcal B$-saturated or not.
For $M = S(2)$, we have $\Ext^1(B,M) \neq 0$ for any bristle $B$. For $M = B_c$,
we have $\Ext^1(B_c,M) \neq 0$. This shows that the bristled modules
$S(2)$ and $B_c$ are not $\mathcal B$-saturated. On the other hand, we have
$\Ext^1(X,Y) = 0$ for $X$ regular and $Y$ preinjective.
This implies that the preinjective modules are $\mathcal B$-saturated modules.
It follows that the only $\mathcal B$-saturated bristled modules are the modules
$M = I_t$ with $0\le t \le q$.
\end{proof}

\subsection{} 
{\bf Remark.} There is the following consequence:
{\it For $n = 2$, the set of indecomposable bristled modules is characterized
by their dimension vectors.}

Let us mention already here that for $n\ge 3$, the set of indecomposable bristled modules
(even the set of faithful ones)
is not characterized by their dimension vectors, as already the dimension vector $(3,2)$ shows:
$$
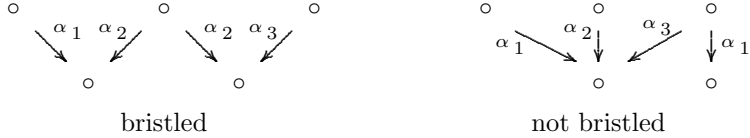
\hbox{\beginpicture
  \setcoordinatesystem units <1cm,1cm>
  \put{\beginpicture
  \multiput{$\circ$} at 0 1  1 0  2 1  3 0  4 1 /
  \arr{0.3 0.7}{0.7 0.3}
  \arr{1.7 0.7}{1.3 0.3}
  \arr{2.3 0.7}{2.7 0.3}
  \arr{3.7 0.7}{3.3 0.3}
  \put{$\ssize\alpha_1$} at 0.7 0.7
  \put{$\ssize\alpha_2$} at 1.3 0.7
  \put{$\ssize\alpha_2$} at 2.7 0.7
  \put{$\ssize\alpha_3$} at 3.3 0.7
  \put{bristled} at 2 -.5
  \endpicture} at 0 0
  \put{\beginpicture
  \multiput{$\circ$} at 0 1  1.5 0  1.5 1  3 0  3 1 /
  \arr{0.4 0.7}{1.2 0.3}
  \arr{1.5 0.7}{1.5 0.3}
  \arr{2.6 0.7}{1.9 0.3}
  \arr{3 0.7}{3 0.3}
  \put{$\ssize\alpha_1$} at 0.3 0.5
  \put{$\ssize\alpha_2$} at 1.2 0.7
  \put{$\ssize\alpha_3$} at 2.25 0.75
  \put{$\ssize\alpha_1$} at 3.3 0.5
  \put{not bristled} at 1.5 -.5
  \endpicture} at 6 0
\put{} at 0 -1.1
  \endpicture}
  $$

\section{Proof of the main theorem.}\label{proof}

In this section, we will assume that $n\ge 3$. We will identify
the index set $\{1,2,\dots,n\}$ of the arrow set of $K(n)$ with $\mathbb Z/n\mathbb Z$, thus
sometimes we will write $\alpha_{n+1}$ instead of $\alpha_1$.  

\subsection{}
We will use the following notation for some bristles (this convention
deviates from the
previous notation, but there should be no confusion).
If $1\le r \le n$, let $B(r)$ be defined by $B(r)_1 = B(r)_2 = k$
with $\alpha_i$ being the identity map in case $i = r$, and the zero map otherwise.
If $r\neq s$ are numbers between $1$ and $n$, let $B(r,s)$ be defined by
$B(r,s)_1 = B(r,s)_2 = k$
with $\alpha_i$ being the identity map in case $i \in\{r,s\}$, and the zero map otherwise.

\subsection{}\label{choice}
Let $\mathcal B_0$ be the set of (the isomorphism classes of)
the bristles $B(n\!-\!1), B(n)$ and the 
bristles $B(r,r+1)$ with $1\le r \le n$ (recall that $B(n,n+1)$ is just $B(n,1)$).
Note that the cardinality of $\mathcal B_0$ is $n+2$.
Let $\mathcal B'_0$ be obtained from $\mathcal B_0$ by deleting $B(n-1,n)$, thus 
the cardinality of $\mathcal B'_0$ is $n+1$.

\subsection{}
For the proof of the assertion (b) of Theorem \ref{theorem},
we may assume that $M$ is indecomposable. In
     case $M$ is preprojective, nothing has to be shown since $\tau^tM = 0$ for $t$ large.
     The case that $M$ is preinjective is considered in assertion (a).
     Finally, there is the case that $M$ is regular.
     Since the class of modules which are generated
     by $\mathcal B_0$ and $\mathcal B$-saturated is closed under extensions, it is sufficient
     to consider just the elementary modules as introduced by Crawley-Boevey (see
     \cite{[K1]} and \cite{[K2]}, or section \ref{elem-gen} below). 
   The proof of (a) dealing with the preinjective modules
     and of the statement (b) for elementary modules proceeds in
     several steps and will be completed in \ref{endofproof}.

\subsection{} {\it The preinjective module $I_2$ is generated by $\mathcal B'_0$, thus by
$\mathcal B_0$.}\label{I-two}
	\medskip

For the proof of \ref{I-two} (and also for the proof of \ref{reg})
we are going to use the universal cover $\widetilde K(n)$ of the $n$-Kronecker quiver 
$K(n)$
(see for example \cite{[FR1]}
and \cite{[FR2]}) in order to exhibit sufficiently many bristles as submodules
of the modules $I_2$ and $\tau B(1)$.
Note that the quiver $\widetilde K(n)$ is
the $n$-regular tree with bipartite orientation. The representations of $\widetilde K(n)$
will be called $\widetilde K(n)$-modules (we may consider them as modules over the
corresponding path algebra; we should note that the path algebra of such a quiver
with infinitely many vertices does not have a unit element,
but at least sufficiently many primitive idempotents).

The quiver $\widetilde K(n)$ is a ``labeled'' quiver:
The
projection $\pi\colon \widetilde K(n) \to K(n)$ sends the arrow $\alpha$ of $\widetilde K(n)$ to an arrow $\pi(\alpha)$ of $K(n)$, 
this is one of $\alpha_1,\dots,\alpha_n,$ and we use $\pi(\alpha)$ (or its index $i$
with $1\le i \le n$)
as a label for $\alpha$. By abuse of notation, we sometimes denote $\alpha$ by
its label $\pi(\alpha)$.
We also denote by $\pi$ the corresponding push-down functor $\mo \widetilde K(n) \to
\mo K(n)$.

\begin{proof}[Proof of \rm \ref{I-two}]

Let $M = I_2$ and let 
$X$ be a representation of $\widetilde K(n)$ with 
$\pi(X) = M,$ thus there is a vertex $z$ of $\widetilde K(n)$
such that the support $Q$ of $X$ consists of the vertices in $\widetilde K(n)$ 
with distance at most $2$
from $z$ and $\dim X_z = n-1$, whereas $\dim X_a = 1$ 
for the remaining vertices $a$ in $Q$.

Here is, in the case $n=4,$ the quiver $Q$ as well as the dimension vector of $X$
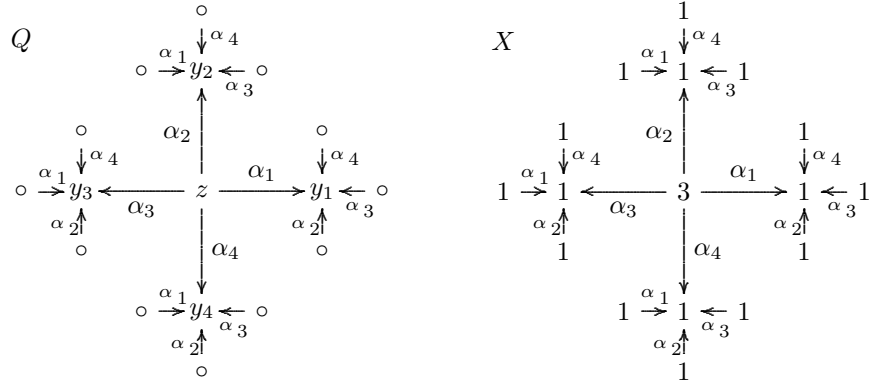
$$
\hbox{\beginpicture
  \setcoordinatesystem units <.8cm,.8cm>
  \put{\beginpicture
  \put{$z$} at 0 0
  \put{$y_1$} at 2 0 
\put{$y_2$} at 0 2
\put{$y_3$} at -2 0
\put{$y_4$} at 0 -2
\multiput{$\circ$} at 
    1 2  0 3  -1 2
    2 1  2 -1  3 0
    -2 1  -2 -1  -3 0
    1 -2  0 -3  -1 -2 /
\arr{0.3 0}{1.7 0}
\arr{-.3 0}{-1.7 0}
\arr{0 0.3}{0 1.7}
\arr{0 -.3}{0 -1.7}

\arr{2.7 0}{2.3 0}
\arr{2 0.7}{2 0.3}
\arr{2 -.7}{2 -.3}

\arr{-2.7 0}{-2.3 0}
\arr{-2 0.7}{-2 0.3}
\arr{-2 -.7}{-2 -.3}

\arr{0 2.7}{0 2.3}
\arr{0.7 2}{0.3 2}
\arr{-.7 2}{-.3 2}

\arr{0 -2.7}{0 -2.3}
\arr{0.7 -2}{0.3 -2}
\arr{-.7 -2}{-.3 -2}
\put{$ \alpha_1$} at 1 0.3
\put{$ \alpha_2$} at -.4 .95
\put{$ \alpha_3$} at -1 -.3
\put{$ \alpha_4$} at .4 -1

\multiput{$\ssize \alpha_1$} at -.5 2.25  -2.5 0.25  -.5 -1.75 /
\multiput{$\ssize \alpha_2$} at 1.7 -.6 -.3 -2.6  -2.3  -.6 /
\multiput{$\ssize \alpha_3$} at 2.6 -.3  0.6 1.7  0.5 -2.3 /
\multiput{$\ssize \alpha_4$} at .35 2.55  -1.65 0.55  2.35 0.55 /
\put{$Q$} at -3 2.5
\endpicture} at 0 0
\put{\beginpicture
\put{$3$} at 0 0
\multiput{$1$} at 2 0  0 2  -2 0  0 -2
    1 2  0 3  -1 2
    2 1  2 -1  3 0
    -2 1  -2 -1  -3 0
    1 -2  0 -3  -1 -2 /
\arr{0.3 0}{1.7 0}
\arr{-.3 0}{-1.7 0}
\arr{0 0.3}{0 1.7}
\arr{0 -.3}{0 -1.7}

\arr{2.7 0}{2.3 0}
\arr{2 0.7}{2 0.3}
\arr{2 -.7}{2 -.3}

\arr{-2.7 0}{-2.3 0}
\arr{-2 0.7}{-2 0.3}
\arr{-2 -.7}{-2 -.3}

\arr{0 2.7}{0 2.3}
\arr{0.7 2}{0.3 2}
\arr{-.7 2}{-.3 2}

\arr{0 -2.7}{0 -2.3}
\arr{0.7 -2}{0.3 -2}
\arr{-.7 -2}{-.3 -2}
\put{$ \alpha_1$} at 1 0.3
\put{$ \alpha_2$} at -.4 .95
\put{$ \alpha_3$} at -1 -.3
\put{$ \alpha_4$} at .4 -1

\multiput{$\ssize \alpha_1$} at -.5 2.25  -2.5 0.25  -.5 -1.75 /
\multiput{$\ssize \alpha_2$} at 1.7 -.6 -.3 -2.6  -2.3  -.6 /
\multiput{$\ssize \alpha_3$} at 2.6 -.3  0.6 1.7  0.5 -2.3 /
\multiput{$\ssize \alpha_4$} at .35 2.55  -1.65 0.55  2.35 0.55 /
\put{$X$} at -3 2.5
\endpicture} at 8 0
\endpicture}
$$
Always, we denote by $y_j$ the neighbor of $z$ such that
the arrow $z \to y_j$ has label $j$.
Recall that a {\it leaf} of a tree is a vertex with at most one neighbor. The
leaves of $Q$ are just the vertices which have distance $2$ from the center $z$;
in the picture on the left, they are indicated by small circles $\circ$.
If $x$ is a leaf of $Q$ and the arrow starting in $x$ has label $\alpha_i$ and ends in $y_j$,
we write $x = x(j,i)$ and we say that $x$ is {\it of type $i$}

In order to find bristle submodules in $M$, we will exhibit some
submodules of $X$; all these submodules are indecomposable
and their support is a path in $Q$.
          \medskip

For $1\le j \le n$, let $Q(j)$ be the full subquiver of $Q$ with vertices $y_j$ and $x(j,i)$
where $1\le i \le n$ and $i\neq j$. Note that
$Q(j)$ is a quiver with one sink and $n-1$ sources (thus it is a so-called subspace quiver).

For $1\le j \le n$, let $Y(j) = X|Q(j)$. Since the subquiver $Q(j)$ of $Q$ is closed
under successors, we see that $Y(j)$ is a submodule of $X$. Since the subquivers $Q(j)$
are pairwise disjoint, the submodule $Y(j)$ form a direct sum $Y = \bigoplus_j Y(j)$
and the factor module
$X/Y$ is nothing else than $X_z$, considered as a direct sum of copies
of the simple representation $S(z)$.

Let $N(j) = \pi(Y(j))$ and $N = \bigoplus_j N(j)$; these are submodules of $M$.
The factor module $M/N$ is a direct sum of copies of $S(1)$ and may be identifies with $M_z$.
We are going to look on the one hand at $N$, and will show that any $N(j)$ is generated by
bristles which are denoted by $D(x)$ and $E(x)$, where $x$ is a leaf in $Q(j)$, see (I) below. 
On the other hand, we have to find bristle submodules $M(i,j)$ of $M$ which are
not contained in $N$, see (II) below. 
    \medskip

For any leaf $x$, let $P(x)$ be the indecomposable projective representation of $Q$
corresponding to the vertex $x$. Note that $P(x)$ is a bristle, and, of course, also
$D(x) = \pi(P(x))$ is a bristle. 
If $x$ is of type $i$, the bristle $D(x)$ is isomorphic to $B(i)$.
   \medskip

{\bf (I) The bristle submodules $D(x)$ and $E(x)$ of $M$.} We are going to show the following:
	\medskip

{\it For any $j$, the submodule $N(j)$ is generated by bristle submodules
which belong to $\mathcal B'_0$.}
	\medskip 

Proof.
Fix some $j$ and write $[i] = x(j,i).$ In order to reduce the number of brackets, we will
write $P[i]$ instead of $P(x(j,i))$ and so on.

Let $G(j)$ be the set of leaves in $Q(j)$ of types $n-1$ and $n$.
If $x\in G(j)$, then $D(x) = \pi(P(x))$ belongs to $\mathcal B'_0$.
If $j \in \{n-1,n\}$, then $G(j)$ is a single leave, otherwise $G(j)$ consists of two
leaves.

We consider in $Q(j)$ subquivers with three vertices, namely those of the following kind:
$$
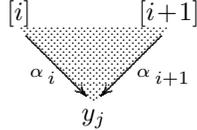
\hbox{\beginpicture
  \setcoordinatesystem units <1cm,1cm>
\put{$[i]$} at 0 1.2
\put{$y_j$} at 1 -.2  
\put{$[i\!+\!1]$} at 2 1.2
\arr{0.1 0.9}{0.9 0.1}
\arr{1.9 0.9}{1.1 0.1}
\put{$\ssize\alpha_i$\strut} at 0.3 0.4
\put{$\ssize\alpha_{i+1}$\strut} at 1.9 0.4
\setshadegrid span <.4mm>
\hshade 0 1 1 <z,z,z,z> 1 0 2 /
\endpicture}
$$
such that $i\neq n\!-\!1$
(note that we also have that $i\neq j$ and $i+1 \neq j$, thus $i\neq j-1$)
and we denote by $V[i] = V(x(j,i))$ the indecomposable
$kQ$-module with this support $\{y_j,[i],[i+1]\}.$
Let $H(j)$ be the set of all $1\le i \le n$
with $i\notin\{j-1,j,n-1\}$. Note that $i\in H(j)$ if and only if $V[i]$
is defined. We claim that

\begin{equation}  \label{eq1}
 Y(j) = \sum_{i\in G(j)} P[i] + \sum_{i\in H(j)} V[i].
\end{equation}

Any $V[i]$ has two bristle submodules, namely $P[i]$ and $P[i+1]$.
The Kronecker module $\pi(V[i])$ has the bristle submodules $D[i] = \pi(P[i])$ and
$D[i+1] = \pi(P[i+1])$, but actually
many more bristle submodules: any non-zero element in the top
of $\pi(V[i])$ generates a bristle.
We are interested in the submodule $E[i] \subset \pi(V[i])$
which is isomorphic to $B(i,i+1)$, thus belongs to $\mathcal B'_0$ (since $i\neq n-1$).
We claim that

\begin{equation} \label{eq2}
 N(j) = \pi(Y(j)) = \sum_{i\in G(j)} D[i] + \sum_{i\in H(j)} E[i].
\end{equation}

For the proof of (\ref{eq1}) and (\ref{eq2})
we distinguish the cases whether $j$ is one of the numbers $n-1, n$ or not.

Consider first the case $j\notin\{n-1,n\}$. Then $G(j) = \{n-1,n\}$ and
$H(j) =\{j+1,j+2,\dots,n-2,n,1,\dots,j-2\}$. The support of the modules $V[i]$
with $i\in H(j)$ are the various shaded triangles in the following picture:
$$
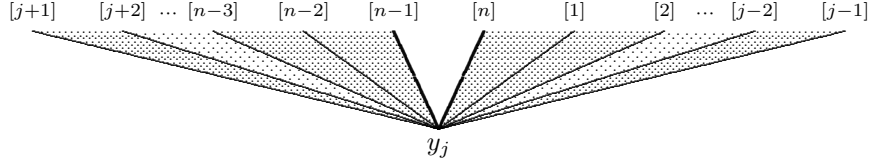
\hbox{\beginpicture
  \setcoordinatesystem units <.6cm,1.3cm>
\put{} at 0 1.4

\put{$\ssize{[j+1]}$} at -9 1.2
\put{$\ssize{[j+2]}$} at -7 1.2
\put{$\ssize{\cdots}$} at -6 1.17

\put{$\ssize{[n-3]}$} at -5 1.2
\put{$\ssize{[n-2]}$} at -3 1.2
\put{$\ssize{[n-1]}$} at -1 1.2
\put{$\ssize{[n]}$} at 1 1.2
\put{$\ssize{[1]}$} at 3 1.2
\put{$\ssize{[2]}$} at 5 1.2
\put{$\ssize{\cdots}$} at 5.9 1.17
\put{$\ssize{[j-2]}$} at 7 1.2
\put{$\ssize{[j-1]}$} at 9 1.2
\put{$y_j$} at 0 -0.2
\plot -9 1  0 0  9 1 /
\plot -7 1  0 0  7 1 /
\plot -5 1  0 0  5 1 /
\plot -3 1  0 0  3 1 /
\plot -.98 1  .02 0  1.02 1 /
\plot -1 1  0 0  1 1 /
\plot -1.02 1  -.02 0  .98 1 /

\setshadegrid span <.4mm>
\hshade 0 0 0 <z,z,z,z> 1 -9 -7 /
\hshade 0 0 0 <z,z,z,z> 1 -5 -1 /
\hshade 0 0 0 <z,z,z,z> 1  1 5 /
\hshade 0 0 0 <z,z,z,z> 1  7 9 /
\setshadegrid span <.6mm>
\hshade 0 0 0 <z,z,z,z> 1 -7 -5 /
\hshade 0 0 0 <z,z,z,z> 1  5 7 /

\endpicture}
$$
Let us show (\ref{eq1}) and (\ref{eq2}) in this case.
Consider first (\ref{eq1}). Of course, the submodules mentioned on the right side of (\ref{eq1})
have support in $Q(j)$, thus they are submodules of $Y(j)$. What we have
to show is that for any leaf $x = [i]$ of $Q(j)$, the module $P(x)$ is included
in the right side. We only have to deal with the vertices $[i]$, where $1\le i \le n-2$
(and, of course,  $i\neq j$).
First, let $1\le i \le j-1$. We know that $P[i]$ is contained in $V[i-1]$ and
$i-1$ is either equal to $n$ or $1\le i-1 \le j-2$, thus always $i-1$ belongs to $H(j)$.
Second, let $j+1 \le i \le n-2$. In this case, we use that  $P[i]$ is contained in $V[i]$
and that $i$ itself belongs to $H(j)$. This completes the proof of (\ref{eq1}).

The equality (\ref{eq1}) implies that
$$
 \pi(Y(j)) = \sum_{i\in G(j)} \pi(P[i]) + \sum_{i\in H(j)} \pi(V[i]),
$$
 and $\pi(P[i]) =  D[i]$. It remains to be seen that for any $i\in H(j)$, the submodule
 $\pi(V[i])$ is contained in the right hand side $R$ of (\ref{eq2}). For $i = n$, we have
 $\pi(V[n]) = D[n] + E[n]$ and $n$ belongs both to $G(j)$ and to $H(j)$.
 Next, let $1\le i \le j-2.$ We use induction on $i$ in order to see that $\pi(V[i])$
 is contained in $R$. Now
$$
  \pi(V[i]) = \pi(P[i])+ E[i]\subset \pi(V[i-1]) + E[i].
$$
  For $i = 1$, we know already that
  $\pi(V[i-1]) = \pi (V[n])$ belongs to $R$; for
  $i > 1$, we know by induction that $\pi V[i-1] \subset R$. Thus, we see that
  let all $1\le i \le j-2,$ the module $\pi(V[i])$ is contained in $R$.

Second, we have to show that $\pi(V[i])$ is contained in $R$, for $j+1 \le i \le n-2.$
This time, we use induction going downwards, starting with $i = n-2$,
and the fact that $\pi(V[i]) = E[i]+\pi(P[i+1]).$
For the start $i= n-2$ of the induction, we have $i+1 = n-1$ and this number belongs to
$G(j)$, thus $\pi(P[i+1]) = D[i+1] = D[n-1]$ belongs to $R$.
For $i < n-2$, we have $\pi(P[i+1]) \subset \pi(V[i+1]$
and $\pi(V[i+1] \subseteq R$, by induction. This completes the proof of
(\ref{eq2}) in the case that $j \notin\{n-1,n\}.$
	    \medskip

Now, we deal with the two cases where $j\in\{n-1,n\}$. In both cases, the proof
will be similar to the previous proof, thus we only describe the setting in detail.

If $j = n-1,$ then $G(j) = \{n\}$ and
$H(j) =\{n,1,\dots,n-3\}$. The support of the modules $V[i]$
with $i\in H(j)$ are the various shaded triangles in the following picture:
$$
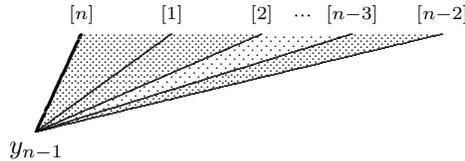
\hbox{\beginpicture
  \setcoordinatesystem units <.6cm,1.3cm>
\put{$\ssize{[n]}$} at 1 1.2
\put{$\ssize{[1]}$} at 3 1.2
\put{$\ssize{[2]}$} at 5 1.2
\put{$\ssize{\cdots}$} at 5.9 1.17
\put{$\ssize{[n-3]}$} at 7 1.2
\put{$\ssize{[n-2]}$} at 9 1.2
\put{$y_{n-1}$} at 0 -0.2
\plot 0 0  9 1 /
\plot 0 0  7 1 /
\plot 0 0  5 1 /
\plot 0 0  3 1 /
\plot .02 0  1.02 1 /
\plot  0 0  1 1 /
\plot  -.02 0  .98 1 /

\setshadegrid span <.4mm>
\hshade 0 0 0 <z,z,z,z> 1  1 5 /
\hshade 0 0 0 <z,z,z,z> 1  7 9 /
\setshadegrid span <.6mm>
\hshade 0 0 0 <z,z,z,z> 1  5 7 /

\endpicture}
$$

For $j = n,$ we have $G(j) = \{n-1\}$ and
$H(j) =\{1,\dots,n-2\}$. The support of the modules $V[i]$
with $i\in H(j)$ are the various shaded triangles in the following picture:
$$
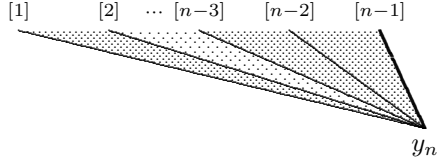
\hbox{\beginpicture
  \setcoordinatesystem units <.6cm,1.3cm>
\put{$\ssize{[1]}$} at -9 1.2
\put{$\ssize{[2]}$} at -7 1.2
\put{$\ssize{\cdots}$} at -6 1.17

\put{$\ssize{[n-3]}$} at -5 1.2
\put{$\ssize{[n-2]}$} at -3 1.2
\put{$\ssize{[n-1]}$} at -1 1.2
\put{$y_n$} at 0 -0.2
\plot -9 1  0 0  /
\plot -7 1  0 0  /
\plot -5 1  0 0  /
\plot -3 1  0 0  /
\plot -.98 1  .02 0  /
\plot -1 1  0 0  /
\plot -1.02 1  -.02 0  /

\setshadegrid span <.4mm>
\hshade 0 0 0 <z,z,z,z> 1 -9 -7 /
\hshade 0 0 0 <z,z,z,z> 1 -5 -1 /
\setshadegrid span <.6mm>
\hshade 0 0 0 <z,z,z,z> 1 -7 -5 /

\endpicture}
$$

As we have mentioned, in both case $j = n-1$ and $j = n$, the proof
of the equalities (\ref{eq1}) and (\ref{eq2}) proceeds as in the case
when $j \notin\{n-1,n\}$, but using just one induction and not two.
In the case $j = n-1$, we have to use the upgoing induction, in the case $j = n$
the downgoing induction.
    \bigskip

{\bf (II) The bristle submodules $M(i,j)$ of $M$.} 
For any pair $i\neq j$ with $1\le i,j\le n$,
there is a (unique) path $p(i,j) = \alpha_i^{-1}\alpha_j\alpha_i^{-1}\alpha_j$ in $Q$
between $x(j,i)$ and $x(i,j)$:
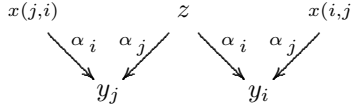
$$
\hbox{\beginpicture
  \setcoordinatesystem units <1cm,1cm>
\put{} at 0 -0.3
\put{$\ssize {x(j,i)}$} at 0 1.1
\put{$y_j$} at 1 0  
\put{$z$} at 2 1.1
\put{$y_{i}$} at 3 0  
\put{$\ssize {x(i,j)}$} at 4 1.1
\arr{0.2 0.8}{0.8 0.2}
\arr{1.8 0.8}{1.2 0.2}
\arr{2.2 0.8}{2.8 0.2}
\arr{3.8 0.8}{3.2 0.2}
\multiput{$\ssize\alpha_{i}$\strut} at 0.65 0.7  2.65 0.7 /
\multiput{$\ssize\alpha_j$\strut} at 1.3 0.7  3.3 0.7 /
\endpicture}
$$
Let $W(i,j)$ be the indecomposable $kQ$-module whose support is
the path $p(i,j)$.

Observe that $W(i,j)$ is a submodule of $X$; it is defined by
$$
 W(i,j)_z = \bigcap_{s\, \notin\{i,j\}} \Ker(\alpha_s)
$$
(this is a one-dimensional
subspace of $X_z$) and $W(i,j)_a = X_a$ for the remaining vertices $a$ on the path $p(i,j)$.

Let us look at $\pi(W(i,j))$.
We claim that the $n$-Kronecker module $\pi(W(i,j))$ 
contains a bristle $M(i,j)$ which is isomorphic to 
$B(i,j)$. The top of $\pi(W(i,j))$ is given by
$$
  W(i,j)_{x(j,i)}\oplus W(i,j)_z\oplus W(i,j)_{x(i,j)} =
  k\oplus k\oplus k.
$$
The element $(1,1,1)\in k^3$ is mapped under 
both $\alpha_{i}$ and $\alpha_{j}$ to
$(1,1)\in W(i,j)_{y_{i}}\oplus W(i,j)_{y_j}$. This shows that 
the $kK(n)$-submodule $M(i,j)$ of $\pi(W(i,j))$
generated by $(1,1,1) \in 
W(i,j)_{x(j,i)}\oplus W(i,j)_z\oplus W(i,j)_{x(i,j)}$ is isomorphic to $B(i,j)$.
Altogether, we see that 
$$
 \pi(W(i,j)) = \pi(P(x(j,i))) + M(i,j)+\pi(P(x(i,j))).
$$
Since both $\pi(P(x(j,i)))$ and $\pi(P(x(i,j)))$ are contained in $N$, it
follows that 
\begin{equation} \label{eq3}
 \pi(W(i,j)) \ \subset\ N + M(i,j).
\end{equation}
	\medskip

The pairs $(i,j)$ which we will consider are those of the form $(i,i+1)$. 
Let $I$ be a subset of $\{1,2,\dots,n\}$ of cardinality $n-1$, thus obtained from
$\{1,2,\dots,n\}$ by deleting one of its elements. We will use the submodules $W(i,i+1)$ 
with $i\in I$.

For example, for $n=4$, and $I = \{1,2,3\}$, here are the dimension vectors 
of the submodules $W(i,i+1)$ of $X$ with $i\in I$
(as well as the corresponding paths 
$\alpha_i{}^{-1}\alpha_{i+1}\alpha_i{}^{-1}\alpha_{i+1}$):
$$
\hbox{\beginpicture 
  \setcoordinatesystem units <.65cm,.65cm>
\put{} at 0 4
\put{\beginpicture
\multiput{$0$} at  -2 0 
    1 2  0 3  
     3 0
    -2 1  -2 -1  -3 0
    1 -2  0 -3
        2 1  0 -2  -1 -2 /
\multiput{$1$} at 0 0  2 0
    -1 2  0 2  2 -1 /
\setdots <1mm>
\arr{-.3 0}{-1.7 0}
\arr{0 0.3}{0 1.7}

\arr{2.7 0}{2.3 0}
\arr{2 -.7}{2 -.3}

\arr{-2.7 0}{-2.3 0}
\arr{-2 0.7}{-2 0.3}
\arr{-2 -.7}{-2 -.3}

\arr{0 2.7}{0 2.3}
\arr{0.7 2}{0.3 2}
\arr{-.7 2}{-.3 2}

\arr{0 -2.7}{0 -2.3}
\arr{0.7 -2}{0.3 -2}

\arr{0 -.3}{0 -1.7}
\arr{-.7 -2}{-.3 -2}

\put{$ \alpha_1$} at 1 0.3
\put{$ \alpha_2$} at -.4 1

\setsolid
\arr{0.3 0}{1.7 0}
\arr{0 .3}{0 1.7}
\arr{2 -.7}{2 -.3}
\arr{-.7 2}{-.3 2}

\multiput{$\ssize \alpha_1$} at  -.5 2.4 /
\multiput{$\ssize \alpha_2$} at   1.6 -.6 /

\put{$\alpha_1{}^{-1}\alpha_2\alpha_1{}^{-1}\alpha_2$} at 0 -4
\put{$W(1,2)$} at -2.4 3
\endpicture} at 0 0

\put{\beginpicture
\multiput{$0$} at  2 0  2 -1  -1 2 
    2 1  3 0
    -2 1  -3 0
    1 -2 -1 -2
        0 3  0 -2  0 -3 /
\multiput{$1$} at 0 0  0 2
     -2 0   -2 -1  1 2 /

\setdots <1mm>
\arr{0.3 0}{1.7 0}

\arr{2.7 0}{2.3 0}
\arr{2 0.7}{2 0.3}
\arr{2 -.7}{2 -.3}
\arr{-2 0.7}{-2 0.3}
\arr{-2.7 0}{-2.3 0}

\arr{-.7 2}{-.3 2}
\arr{0.7 -2}{0.3 -2}

\arr{-.7 -2}{-.3 -2}
\arr{0 -2.7}{0 -2.3}
\arr{0 2.7}{0 2.3}
\arr{0 -.3}{0 -1.7}

\setsolid
\arr{0 0.3}{0 1.7}
\arr{-.3 0}{-1.7 0}
\arr{.7 2}{.3 2}
\arr{-2 -.7}{-2 -.3}

\put{$ \alpha_2$} at -.4 .95
\multiput{$\ssize \alpha_2$} at -2.4 -.55  /
\put{$ \alpha_3$} at -.9 -.3
\multiput{$\ssize \alpha_3$} at   .5 1.7 /

\put{ $\alpha_2{}^{-1}\alpha_3\alpha_2{}^{-1}\alpha_3$} at 0 -4

\put{$W(2,3)$} at -2.4 3
\endpicture} at 7 0
\put{\beginpicture
\multiput{$0$} at 0 2  2 0  3 0  0 -3  -3 0 
    1 2  0 3  -1 2
    2 1  2 -1 
    -2 -1 
    -1 -2 /
\multiput{$1$} at 0 0  -2 0  0 -2  1 -2  -2 1 /
\setdots <1mm>
\arr{0 0.3}{0 1.7}

\arr{2 0.7}{2 0.3}
\arr{2 -.7}{2 -.3}

\arr{-2 -.7}{-2 -.3}

\arr{0 2.7}{0 2.3}
\arr{0.7 2}{0.3 2}
\arr{-.7 2}{-.3 2}

\arr{0 -2.7}{0 -2.3}
\arr{-.7 -2}{-.3 -2}
\arr{0.3 0}{1.7 0}
\arr{2.7 0}{2.3 0}
\arr{-2.7 0}{-2.3 0}

\setsolid
\put{$ \alpha_4$} at 0.45 -1
\put{$ \alpha_3$} at -1 -.3
\arr{-2 0.7}{-2 0.3}
\arr{0.7 -2}{0.3 -2}
\arr{0 -.3}{0 -1.7}
\arr{-.3 0}{-1.7 0}

\multiput{$\ssize \alpha_4$} at   -1.6 .5 /
\multiput{$\ssize \alpha_3$} at 0.5 -2.3   /

\put{$\alpha_3{}^{-1}\alpha_4\alpha_3{}^{-1}\alpha_4$} at 0 -4
\put{$W(3,4)$} at -2.4 3
\endpicture} at 14 0
\endpicture}
$$
	\medskip 

Here is the essential assertion:
{\it If $I \subset\{1,2,\dots,n\}$ has cardinality $n-1$, then}
\begin{equation}  \label{eq-final}
        M = N+\sum_{i\in I} M(i,i+1).
\end{equation}

\begin{proof}[Proof of \rm{(\ref{eq-final})}] 
We can assume that $I = \{1,2,\dots,n-1\}$. First, we show that 
there is a vector space decomposition
\begin{equation}  \label{eq4}
        X_z = \bigoplus_{i\in I} W(i,i+1)_z. 
\end{equation}

For the proof of (\ref{eq4}), we recall 
that the arrow of type $\alpha_i$ starting in $z$ 
ends in the vertex $y_i$.
The Auslander-Reiten sequence
in $\mo kQ$ ending in $S(z)$ is of the form
$$
\hbox{\beginpicture
  \setcoordinatesystem units <1cm,.4cm>
\put{$0\strut$} at 0.5 0
\put{$X$\strut} at 1.5 0  
\put{$\bigoplus_{1\le i \le n} I(y_i)\strut$} at 4 0  
\put{$S(z)\strut$} at 8 0  
\put{$0$\strut} at 9.4 0
\arr{0.7 0}{1.3 0}
\arr{1.7 0}{2.8 0}
\arr{5.2 0}{7.3 0}
\arr{8.5 0}{9.2 0}
\put{$\iota$} at 2.25 0.5  
\put{$\ssize{[\pi(1),\dots,\pi(n)]}$} at 6.2 0.5  
\endpicture}
$$
where the support of $I(y_i)$ consists of $y_i$ and its neighbors, where
$\pi(i)\colon  I(y_i) \to S(z)$ is the canonical projection. We
can assume that $\iota$ is an inclusion map. At the vertex $z$, there is the
corresponding exact sequence of vector spaces
$$
\hbox{\beginpicture
  \setcoordinatesystem units <1cm,.4cm>
\put{$0$\strut} at 0.5 0
\put{$X_z$\strut} at 1.5 0  
\arr{0.7 0}{1.2 0}
\arr{1.8 0}{2.8 0}
\put{$\iota_z$} at 2.25 0.5

\put{$\bigoplus_{1\le i \le n} I(y_i)_z$\strut} at 4 0  
\put{$S(z)_z$\strut} at 9 0  
\put{$0$\strut} at 10.4 0
\arr{5.2 0}{8.3 0}
\arr{9.5 0}{10.2 0}
\put{$\ssize{[\pi(1)_z,\dots,\pi(n)_z]}$} at 6.7 0.5  
\endpicture}
$$
with $I(y_i)_z = k,$ $S(z) = k$ and $\pi(i)_z$ the identity map. Thus
             $X_z$ is the kernel of the map
$$
  [\pi(1)_z,\dots,\pi(n)_z] = [1,\dots,1]\colon k^n \to k
$$
and therefore generated by the elements $e(i)-e(i+1)$, with $1\le i < n$; here
$e(i)$ is the canonical generator of $k = I(y_i)_z.$ Of course,
$e(i)-e(i+1)$ belongs to the kernel of $\alpha_s\colon X_z \to X_{a_s}$ for $s\notin\{i,j\}$,
thus to $W(i,i+1)_z$. This shows that
$$
  X_z = \bigoplus_{1\le i < n}\langle e(i)-e(i+1)\rangle
      = \bigoplus_{1\le i < n} W(i,i+1)_z.
$$
This completes the proof of (\ref{eq4}).
        \medskip

As a consequence of (\ref{eq4}), we have 
$$
   X = Y + \sum_{i\in I} W(i,i+1).
$$
We apply $\pi$ and obtain 
$$  
  M = \pi(X) = N + \sum_{i\in I} \pi(W(i,i+1)).
$$
Now we use (\ref{eq3}) in order to see that 
$$
  M = N + \sum_{i\in I} \pi(W(i,i+1)) \subseteq N +  \sum_{i\in I} M(i,i+1) \subseteq M,
$$
thus we obtain (\ref{eq-final}).   
\end{proof}

In order to complete the proof of \ref{I-two}, we choose $I = \{1,2,\dots,n-2,n\}$.
Then all the submodules $M(i,i+1)$ with $i\in I$ belong to $\mathcal B'_0$. We
know already by (I) that $N$ is generated by $\mathcal B'_0$. Thus the equality 
(\ref{eq-final}) shows that $M$ is generated by $\mathcal B'_0$. 
\end{proof}

\subsection{}\label{book}{\bf Remark: Book-keeping.}
For the sake of book-keeping, let us mention for any bristle type the number
of bristles used for generating $M = I_2$. The number of leaves of $Q$ of types $n-1$
and $n$ is $n-1$. Thus, as bristles of the form $D(x)$ we have used $n-1$
bristles isomorphic to $B(n\!-\!1)$ as well as $n-1$ bristles isomorphic to $B(n)$.
Let $I = \{1,2,\dots,n-2,n\}$. For any $i\in I$, 
there are $n-2$ modules of the form $V[i]$, namely with
support inside $Q(j)$, where $j\notin\{i,i+1\}$, thus we have used $n-2$ 
bristles of the form $E[i]$ and these bristles are isomorphic to $B(i,i+1)$.
Finally, for $i\in I$, we needed the module $M(i,i+1)$. Now, $M(i,i+1)$ again
is isomorphic to $B(i,i+1)$. Altogether, we see that for
any bristle $B\in \mathcal B'_0$, the number of bristles isomorphic to $B$ used
in order to generate $I_2$ is $n-1$. 

Therefore, the total number of bristles used for generating $I_2$ is $(n+1)(n-1)
= n^2-1$. Since $\bdim I_2 = (n^2-1,n)$, the module $\topp I_2$ is of length $n^2-1$,
thus a minimal number of bristle submodules which generate $I_2$ consists of 
precisely $n^2-1$ bristles. 
We also should mention that for any bristle $B$, 
$$
 \dim \Hom(B,I_2) = \langle \binom 1 1, \binom{n^2-1}n\rangle =  n^2-1 + n - n^2 = n-1
$$
(since $\Ext^1(B,I_2) = 0).$
Looking at any minimal set of bristle submodules which generate $I_2$, this
set can contain at most $n-1$ bristles isomorphic to $B$. 


\subsection{}\label{reg}
{\it The module $\tau B(1)$
is generated by $\mathcal B_0$ and it has a proper submodule isomorphic to $B(1).$}

\begin{proof}
The proof of the first assertion
is similar to the considerations in \ref{I-two} dealing with $I_2$. Actually, one may 
consider $\tau B(1)$ as a submodule of $I_2$, generated by some of the bristles
constructed above.

Here, we consider
a representation $X'$ of $\widetilde K(n)$ with $\pi(X') = \tau B(1).$ We denote by $Q'$ 
its support: it is obtained from $Q$ defined above as follows:
Deleting the arrow $z \to y_1$, we obtain two connected 
components and we delete the component which contains $y_1$ (this is just $Q(1)$ as considered
in \ref{I-two}). 
Thus $Q'$  consists of the vertices $a$ of $\widetilde K(n)$ 
which have  
distance at most $2$ from $z$, and such that the path 
from $z$ to $a$ does not start with the arrow with label $\alpha_1.$ We have
$\dim X'_z = n-2$, and $\dim X'_a = 1$ for the remaining vertices $a$ in $Q'$.

The module $X'$ can be described as a submodule of $X$ as follows:  
$X'_z$ is the kernel of the  map $X_z \to X_{y_1}$, and $X'_a = X_a$ for the
remaining vertices $a$ in $Q'$. 

As in the proof of  \ref{I-two}, we first consider the subquivers $Q(j)$ and
the restriction of $X'$ to $Q(j)$. Of course, now we deal only with $j\neq 1$.
Actually, the restrictions of $X'$ to $Q(j)$ is just $Y(j)$ as considered in 
\ref{I-two}. We recall that all the modules $Y(j)$ are generated by $\mathcal B_0$.

Second, we need bristles which are not contained in $\pi(Y')$, where $Y'
= \bigoplus_{j=2}^n Y(j)$. We start with the submodules $W(i,i+1) \subset X'$, where 
$2\le i \le n-1$ and take the corresponding submodules $M(i,i+1)$ of $\pi(W(i,i+1))$.
As we know, $M(i,i+1)$ is isomorphic to $B(i,i+1) \in \mathcal B_0$.

Since $\tau B(1) = \pi(Y')+\sum_{i=2}^{n-1}M(i,i+1)$, we see that $\tau B(1)$ is
generated by $\mathcal B_0$. This yields the first assertion.

For the proof of the second assertion, just take any submodule $\pi(P(x(j,1))$ with
$j\neq 1.$ Of course, one may also refer to the Auslander-Reiten formula
which yields non-zero homomorphisms (and this embeddings) $B(1) \to \tau B(1)$, since
obviously $\Ext^1(B(1),B(1))\neq 0.$
\end{proof}
    
Here are, for $n=4$, on the left the dimension vector of $X'$, in the middle and on the right 
those of $W(2,4)$ and $W(3,4)$, respectively. 

$$
\hbox{\beginpicture
  \setcoordinatesystem units <.8cm,.8cm>
\put{\beginpicture
\put{$2$} at 0 0

\multiput{$1$} at  0 2  -2 0  0 -2
    1 2  0 3  -1 2
    -2 1  -2 -1  -3 0
    1 -2  0 -3  -1 -2 /
\arr{-.3 0}{-1.7 0}
\arr{0 0.3}{0 1.7}
\arr{0 -.3}{0 -1.7}

\arr{-2.7 0}{-2.3 0}
\arr{-2 0.7}{-2 0.3}
\arr{-2 -.7}{-2 -.3}

\arr{0 2.7}{0 2.3}
\arr{0.7 2}{0.3 2}
\arr{-.7 2}{-.3 2}

\arr{0 -2.7}{0 -2.3}
\arr{0.7 -2}{0.3 -2}
\arr{-.7 -2}{-.3 -2}

\put{$ \alpha_2$} at -.4 .95
\put{$ \alpha_3$} at -1 -.3
\put{$ \alpha_4$} at .4 -1

\multiput{$\ssize \alpha_1$} at -.5 2.25  -2.5 0.25  -.5 -1.75 /
\multiput{$\ssize \alpha_2$} at -.3 -2.6  -2.3  -.6 /
\multiput{$\ssize \alpha_3$} at 0.6 1.7  0.5 -2.3 /
\multiput{$\ssize \alpha_4$} at .35 2.55  -1.65 0.55  /
\put{$X'$} at -2.7 2.5
\endpicture} at 0 0

\put{\beginpicture
  \setcoordinatesystem units <.55cm,.55cm>
\multiput{$0$} at  -1 2 
    -2 1  -3 0
    1 -2 -1 -2
        0 3  0 -2  0 -3  /
\multiput{$1$} at 0 0
      -2 0  -2 -1  1 2   0 2 /

\setdots <1mm>
\arr{-2 0.7}{-2 0.3}
\arr{-2.7 0}{-2.3 0}

\arr{-.7 2}{-.3 2}
\arr{0.7 -2}{0.3 -2}

\arr{-.7 -2}{-.3 -2}
\arr{0 -2.7}{0 -2.3}
\arr{0 2.7}{0 2.3}
\arr{0 -.3}{0 -1.7}

\setsolid
\arr{0 0.3}{0 1.7}
\arr{-.3 0}{-1.7 0}
\arr{.7 2}{.3 2}
\arr{-2 -.7}{-2 -.3}

\put{$ \alpha_2$} at -.4 .95
\multiput{$\ssize \alpha_2$} at -2.4 -.55  /
\put{$ \alpha_3$} at -.9 -.3
\multiput{$\ssize \alpha_3$} at   .5 1.7 /

\endpicture} at 6 0

\put{\beginpicture
  \setcoordinatesystem units <.55cm,.55cm>
\multiput{$0$} at 0 2  0 -3  -3 0 
    1 2  0 3  -1 2
    -2 -1 
    -1 -2 /
\multiput{$1$} at 0 0  -2 0  0 -2  1 -2  -2 1 /
\setdots <1mm>
\arr{0 0.3}{0 1.7}

\arr{-2 -.7}{-2 -.3}

\arr{0 2.7}{0 2.3}
\arr{0.7 2}{0.3 2}
\arr{-.7 2}{-.3 2}

\arr{0 -2.7}{0 -2.3}
\arr{-.7 -2}{-.3 -2}
\arr{-2.7 0}{-2.3 0}

\setsolid
\put{$ \alpha_4$} at 0.45 -1
\put{$ \alpha_3$} at -1 -.3
\arr{-2 0.7}{-2 0.3}
\arr{0.7 -2}{0.3 -2}
\arr{0 -.3}{0 -1.7}
\arr{-.3 0}{-1.7 0}

\multiput{$\ssize \alpha_4$} at   -1.6 .5 /
\multiput{$\ssize \alpha_3$} at 0.5 -2.3   /
\endpicture} at 11 0
\put{$\alpha_2{}^{-1}\alpha_3\alpha_2{}^{-1}\alpha_3$} at 6 -3
\put{$\alpha_3{}^{-1}\alpha_4\alpha_3{}^{-1}\alpha_4$} at 11 -3
\put{$W(2,3)$} at 4.3 2.5
\put{$W(3,4)$} at 9.3 2.5
\endpicture}
$$
	\medskip

\subsection{} {\bf Remark: Again the book-keeping.} 
The number of leaves of $Q'$ of type $n-1$ is $n-2$, as is the number of leaves 
of type $n$. 
Thus, as bristles of the form $D(x)$ we have used $n-2$
bristles isomorphic to $B(n\!-\!1)$ as well as $n-2$ bristles isomorphic to $B(n)$.
The number of bristles $E(x)$ which we have used are as follows: 
We have used $n-2$ bristles $E(x)$ isomorphic to $B(1,2)$,
we have used $n-3$
bristles $E(x)$ isomorphic to $B(i,i+1)$, for any $i$ with $2\le i \le n-2$,
and we have used $n-2$ bristles $E(x)$ isomorphic to $B(n,1)$.
In addition, there are the bristles $M(i,i+1)$ isomorphic to $B(i,i+1)$. 
Altogether we have used $n-2$ bristles of the form $B(n\!-\!1), B(n)$ and
$B(i,i+1)$ with $1\le i \le n-2$ and with $i=n$ and a single bristle of the
form $B(n-1,n).$

In \ref{reg} we have presented a set of cardinality $n+2$ in order to generate $\tau B(1)$.
As we will see in section \ref{opt2}, there are bristle sets of cardinality $n+1$
which generate $\tau B(1)$, but they are not helpful for the remaining proof of
\ref{theorem}.

\subsection{}
{\bf Corollary to \ref{reg}.} {\it If $B$ is a bristle, then $\tau B$ is generated by bristles, but
is not a bristle,
and $\tau^{-1}B$ is cogenerated by bristles and is not a bristle.}

\begin{proof} Using an automorphism of the Kronecker algebra, we may
    shift $B$ to $B(1)$, apply \ref{reg} and shift back. This shows that $\tau B$
    is generated by bristles. Using duality, it follows that $\tau^{-1}$ is
    cogenerated by bristles.
    \end{proof}
    
\subsection{}\label{ext}
{\it 
If $B,B'$ are non-isomorphic
bristles, then  $\Ext^1(B,\tau B') = 0$.}

\begin{proof} This follows from the isomorphism
$\Ext^1(B,\tau B') \simeq D\Hom(B',B)$,
where $D$ denotes $k$-duality, since
non-isomorphic bristles are $\Hom$-orthogonal.
\end{proof} 

\subsection{}\label{tauM}{\it 
Let $M$ be an indecomposable module and $U$ a submodule of $M$ isomorphic to $B(1)$.
If $\tau(M/U)$ is generated by $\mathcal B_0$, then $\tau M$ is generated by
$\mathcal B_0$.}

\begin{proof} 
We apply $\tau$ to the exact sequence
$$
 0 \to U \to M \to M/U \to 0
$$
and obtain the exact sequence
$$
 0 \to \tau U \to \tau M \to \tau(M/U) \to 0.
$$
Consider any homomorphism $f\colon B \to \tau(M/U)$ with
$B\in \mathcal B_0$. Since $B$ is not isomorphic to $B(1)$, \ref{ext}
shows that $\Ext^1(B,\tau B(1))
= 0$, thus $f$ can be lifted to $\tau M$. And \ref{reg} asserts that
$\tau U$ is generated by $\mathcal B_0$. Since by assumption also
$\tau (M/U)$ is generated by $\mathcal B_0$, it follows
that $\tau M$ is generated by $\mathcal B_0$.
\end{proof} 
     
\subsection{}\label{prein}{\it
All preinjective modules are generated by $\mathcal B_0$
and are $\mathcal B$-saturated.}

\begin{proof}  First, let us show that $I_t$ with $t\ge 0$ is generated by $\mathcal B_0$.
We use induction on $t$.
The assertion is trivial for $t = 0$. The module $I_1$ is obviously generated
by $B(2),\dots,B(n)$ and $B(1,2)$, thus by $\mathcal B_0$. Next, \ref{I-two} asserts 
that $I_2$ is generated by $\mathcal B_0.$ Thus, we can assume that $t\ge 3$. 

Clearly, $\Hom(B(1),I_s) \neq 0$ for any $s\ge 0$, thus $\Hom(B(1),I_{t-2}) \neq 0.$
Take a nonzero 
homomorphism $f\colon B(1) \to I_{t-2}$. Since $t-2 \ge 1$,
the map $f$ has to be a monomorphism
and its cokernel $C$ will be a direct sum of preinjective modules of the form
$I_s$ with $s < t-2.$ We apply $\tau$ to the exact sequence
$$
  0 \to B(1) \to I_{t-2} \to C \to 0
  $$
  and obtain the exact sequence
  $$
    0 \to \tau B(1) \to I_{t} \to \tau C \to 0.
    $$
    Now $\tau C$ is a direct sum of modules of the form $I_{s+2}$ with $s+2 < t.$
By induction, $C$ is generated by $\mathcal B_0$, thus \ref{tauM} shows that
$I_t$ is generated by $\mathcal B_0$.

Of course, all preinjective modules $I_t$ satisfy $\Ext^1(X,I_t) = 0$, for
$X$ regular. Since bristles are regular, we see that $I_t$ is
$\mathcal B$-saturated. 
\end{proof}

\subsection{}\label{homzero}{\it
Let $M$ be an indecomposable  module which is generated by bristles, but
not a bristle. Then $\Hom(M,B) = 0$ for any bristle $B$, thus $\tau M$ is
$\mathcal B$-saturated.}

\begin{proof} Since $M$ is a generated by bristles,
there are bristles $B_1,\dots,B_t$ and a surjective map
$f\colon \bigoplus B_i \to M.$ Assume there exists a non-zero map $g\colon M\to B$.
Then $gf_i \neq 0$ for
some $i$. But a non-zero map between bristles is an isomorphism, thus $g$ is a split monomorphism.
Since $M$ is indecomposable, $g$ has to be an isomorphism, thus $M$ is a bristle. This is the first
assertion. The vector space isomorphism $\Ext^1(B,\tau M) \simeq 
D\underline{\Hom}(M,B)$ shows that $\Ext^1(B,\tau M) = 0$ for any bristle $B$, thus
$\tau M$ is $\mathcal B$-saturated.
\end{proof}

We also will need the dual assertion of \ref{homzero}:
   
\subsection{}\label{dual}{\it 
Let $M$ be an indecomposable module which is cogenerated by bristles, but
not a bristle. Then $\Hom(B,M) = 0$ for any bristle $B$.}

\subsection{}\label{elem-gen}
Let us consider now elementary modules as defined by Crawley-Bovey and studied by
Lukas and Kerner (see for example \cite{[K1]}):
a regular module $M$ is said to be {\it elementary} provided
that for any non-zero regular
submodule $U$ of $M$, the factor module $M/U$ is preinjective. 
	    \medskip
	      
{\it 
Let $M$ be an elementary module with $\Hom(B(1),M) \neq 0.$
Then $\tau^t M$ is generated by $\mathcal B_0$ for all $t\ge 1$.}

\begin{proof}
Let $M$ be an elementary module with $\Hom(B(1),M) \neq 0.$
A non-zero homomorphism $f\colon B(1) \to M$ has to be injective, since $M$ has no
submodule isomorphic to $S(1)$.
Since $M$ is elementary, the cokernel $C$ of $f$ 
has to be preinjective. Thus, we deal with an exact sequence of the form
$$
  0 \to B(1) \to M \to C \to 0.
  $$
  If we apply $\tau$, we obtain the exact sequence
  $$
    0 \to \tau B(1) \to \tau M \to \tau C \to 0.
    $$
Now $\tau C$ is preinjective, thus according to \ref{prein},
generated by $\mathcal B_0$. 
According to \ref{tauM}, also
$\tau M$ is generated by $B_0.$ 

Now, the module  $\tau M$ is again elementary, and it
has $\tau B(1)$ as a submodule. According to \ref{reg}, $\tau B(1)$ has
a submodule isomorphic to $B(1)$. It
follows that $\Hom(B(1),\tau M) \neq 0.$ This shows, that $\tau M$ satisfies
again the assumption, thus also $\tau^2 M$ is generated by $\mathcal B_0$.
Altogether we see that we can use induction.
\end{proof}

Since bristles are
elementary modules, we see:

\subsection{}\label{cor1}{\bf Corollary.} {\it 
If $B$ is a bristle, then the module
$\tau^t B$ is generated by $\mathcal B_0$, for all $t\ge 1$.}

\subsection{}\label{cor2}{\bf Corollary.} {\it 
If $B$ is a bristle, then
the module $\tau^{-t} B$ is cogenerated by bristles, for all $t\ge 1$.}
    \medskip
    
Corollary \ref{cor2} follows from \ref{cor1} by using duality.

\subsection{}\label{bristles}
{\it 
Any $\tau$-orbit contains at most one bristle.
If $B, B'$ are bristles, then we have $\Hom(\tau^tB,B') = 0$
for $t\ge 1.$}

\begin{proof}
Assume $B$ and $\tau^tB$ are bristles, for some $t\ge 1$.
This implies that $\ssize{\binom11}$ is an eigenvector of $\Phi^t$ with
eigenvalue $1$, where
$\Phi$ is the Coxeter transformation for $K(n)$. It is well-known (and easy to see)
that $1$ is not an eigenvalue of $\Phi^t$ for any $t\ge 1$.
This yields the first assertion.

For the proof of the second assertion, let $t\ge 1.$
We know that $\tau^tB$ is generated by bristles, thus there
are maps $f_i\colon B_i\to \tau^tB$ with $B_i$ a bristle, such that $(f_i)_i\colon 
\bigoplus B_i \to \tau^t B$ is surjective. Assume there exists a non-zero
map $g\colon \tau^tB \to B'.$ Then $gf_i \neq 0$ for some $i$. A non-zero map
between bristles is an isomorphism, thus $gf_i$ is an isomorphism and therefore
$g\colon \tau^tB \to B'$ is split epi. Since $\tau^tB$ is indecomposable,
we see that $g$ is an isomorphism, thus $\tau^tB$ is a bristle.
But the $\tau$-orbit of $\tau$ cannot contain a second bristle.
\end{proof}

\subsection{}\label{elem-sat}
{\it 
Let $M$ be an elementary module
    	       with $\Hom(B(1),M) \neq 0.$
Then $\tau^t M$ is $\mathcal B$-saturated for all $t\ge 2.$
}

\begin{proof}
Let $M$ be an elementary module
 with $\Hom(B(1),M) \neq 0.$ Let $t\ge 2$. By \ref{elem-gen} we know that
 $\tau^{t-1}M$ is generated by bristles. If we know that $\tau^{t-1}M$
 is not a bristle, then we can apply \ref{homzero} in order to see that
 $\tau^tM$ is $\mathcal B$-saturated.

Thus assume that $B =\tau^{t-1}M$ is a bristle. Then $M = \tau^{-t+1}B$
is cogenerated by bristles, according to Corollary \ref{cor2}.
Since $\Hom(B(1),M) = 0$, \ref{dual} shows that $M$ has to be a bristle.
But according to the first assertion \ref{bristles}, we know that any $\tau$-orbit
contains at most one bristle.
\end{proof}

\subsection{}\label{cor3}{\bf Corollary.} {\it 
If $B$ is a bristle, then
$\tau^t B$ is $\mathcal B$-saturated for all $t\ge 2.$}
	\medskip

Note that the assumption $t\ge 2$ is essential, since 
$\tau B$ is not $\mathcal B$-saturated: we obviously have $\Ext^1(B,\tau B) \neq 0.$

\subsection{Proof of Theorem \ref{theorem}.}\label{endofproof}
Part (a) has been established in \ref{prein}.
     Thus it remains to deal with part (b). As we have mentioned already, it is sufficient to
     look at a module $M$ which is regular.

First, let us consider the special case of an elementary module $M$.
It has been shown by D.~Baer that for any pair $M',M''$ of regular modules
there is some number $s(M',M'')$ such
that $\Hom(M',\tau^sM'') \neq 0$ for all $s\ge s(M',M'')$ (see for example \cite{[K1]}, 10.7).
Thus $\Hom(B(1),\tau^sM) \neq 0$ for all $s \ge s(B(1),M)$.
According to \ref{elem-gen}
and \ref{elem-sat}, we know that $\tau^sM$ is generated by $\mathcal B_0$
and is $\mathcal B$-saturated for all $s\ge s(B(1),M)+2.$

Now any regular module $M$ has a filtration by elementary modules, say
$$
 0 = M_0 \subset M_1 \subset \cdots \subset M_m = M
 $$
 with all $M_i/M_{i-1}$ for $1\le i \le m$ being elementary.
 Let
$$
  t(M) = \max_i s(B(1),M_i/M_{i-1}) + 2.
$$  Let $t\ge t(M)$.
 The module $\tau^tM$ has the filtration
$$
  0 = \tau^tM_0 \subset \tau^tM_1 \subset \cdots \subset \tau^tM_m = \tau^tM
$$
  and the factors are $\tau^tM_i/\tau^tM_{i-1} = \tau^t(M_i/M_{i-1})$.
  These are elementary modules which are generated by $\mathcal B_0$ and are
  $\mathcal B$-saturated.
  Since the class of modules which are generated by $\mathcal B_0$ and
  $\mathcal B$-saturated is closed under extensions, it follows that also
  $\tau^tM$ is  generated by $\mathcal B_0$ and
  $\mathcal B$-saturated.
  This completes the proof.
  \hfill$\square$

\subsection{}\label{ind}
{\bf Corollary.} {\it 
There is a regular module $X$ which generates all the
preinjective modules, for example $X = \bigoplus_{B\in \mathcal B_0} B$.
There are also indecomposable regular modules which generate all preinjective modules.}
	                              \medskip 

This provides an explicit answer to 
a question raised by Kerner in 1994 (Problem 10.9 in \cite{[K1]}): to exhibit a regular
module $X$ which generates all the preinjective modules. An existence proof for
such a module was given by Kerner in \cite{[K2]}.

\begin{proof}
The first assertion is a direct consequence of \ref{theorem}.

For the second assertion, let $X$ be 
a suitable extensions $0\to B(1) \to
\overline X \to X \to 0$ with $\overline X$ indecomposable (such an extension exists,
since $\mathcal B_0\cup\{B(1)\}$ is a set of pairwise non-isomorphic bristles,
and pairwise non-isomorphic bristles are pairwise
orthogonal,  and $\Ext^1(B,B') \neq 0$ for any pair $B,B'$ of bristles).
Such a module $\overline X$ has dimension vector
$(n+2,n+2)$.
\end{proof}

\subsection{}\label{tame}
{\bf Remark.} Let us stress that the existence of a (finitely generated!)
regular module which generates all preinjective modules is a special feature for
the cases $n\ge 3$: It is a finiteness result in the wild cases which does not hold
in the tame case $n=2$, see Appendix C.

\subsection{}{\bf Addendum to the Main Theorem.}\label{optimal}
{\it Let $\mathcal B'$ be a set of $n+1$ bristles and let $t\ge 3$.
The preinjective module $I_t$ cannot be generated by $\mathcal B'$.}
Thus, the number $n+2$ occurring in the Main Theorem is best possible.

\begin{proof} First, let $t = 3.$
We use the Euler form $\langle-,-\rangle$ on the Grothendieck group $K_0(\mo\Lambda).$
Note that $I_3$ has the dimension vector $\bdim I_3 = (n^3-2n,n^2-1).$
We have
$$
 \langle \binom 1 1,\binom{n^3-2n}{n^2-1}\rangle = (n^3-2n)+(n^2-1)-n(n^2-1) = n^2-n-1.
$$
Since for any bristle $B$, we have $\Ext^1(B,I_3) = 0,$ it follows that
$$
  \dim \Hom(B,I_3) = n^2-n-1.
$$
Let $X$ be the direct sum of $n+1$ pairwise non-isomorphic 
bristles. It follows that there is a right $X$-approximation of the form $f\colon X^t \to I_3$
with $t\le n^2-n-1$. Now 
$$
 \dim\topp X^t = (n+1)t \le (n+1)(n^2-n-1) = n^3-2n-1.
$$
On the other hand,
$$
  \dim \topp I_3 = n^3-2n,
$$
thus $f$ cannot be surjective. This shows that $X$ does not generate $I_3.$

Now consider any $t\ge 3$. Assume that $I_t$ is generated by $\mathcal B'$.
Since $I_t$ generates $I_3,$ it follows that $\mathcal B'$ generates $I_3$.
But we have shown already that this is not possible. 
\end{proof}

\subsection{}\label{opt2} Let us look at optimal ways to generate $\tau B(1)$, using bristle
submodules. 
	\smallskip

{\it Let $\mathcal B'$ be a set of $n+1$ bristles. 
If $\tau B(1)$ is generated by $\mathcal B'$, then $B(1)$ has to belong to $\mathcal B'.$
The set $\mathcal B'_1$ given by $B(1)$ as well as the bristles $B(i,i+1)$ for $1\le i \le n$
generates $\tau B(1)$.}

\begin{proof} We have $\bdim \tau B(1) = (n^2-n-1,n-1)$. If $B$ is a bristle, then
$$
 \langle \binom 1 1, \binom{n^2-n-1}{n-1}\rangle =  
  (n^2-n-1) + (n-1) - n(n-1) = n-2.
$$
Since bristles are bricks and pairwise orthogonal, the Auslander-Reiten formula shows
that $\dim\Ext^1(B(1),B(1)) = 1$ and $\dim\Ext^1(B,B(1))= 0$, provided $B$ is not isomorphic
to $B(1)$. It follows that $\dim \Hom(B(1),\tau B(1)) = n-1$, but 
$\dim\Hom(B,B(1)) = n-2$  in case $B$ is not isomorphic to $B(1).$

Let $\mathcal B'$ be a set of bristles of cardinality $n+1$ and assume
that $\mathcal B'$ does not contain $B(1)$. 
Since $\dim\Hom(B,\tau B(1)) = n-2$ for any $B \in \mathcal B'$, 
there is a right $\mathcal B'$-approximation of the form $f\colon X^t \to \tau B(1)$
with $t\le n-2$. We have 
$\dim\topp X^t = (n+1)t \le (n+1)(n-2) = n^2-n-2$, and $\dim \topp \tau B(1) = n^2-n-1$.
As a consequence, $f$ cannot be surjective. This shows that $\mathcal B'$ cannot generate
$\tau (B(1).$

The proof of the second assertion is similar to the proofs presented in \ref{I-two} and \ref{reg}.
\end{proof}

As we have seen in the proof, the existence of sufficiently many bristles in $\tau B(1)$
which are isomorphic to $B(1)$ is related to the non-vanishing of $\Ext^1(B(1),\tau B(1))$.
Of course, the only non-trivial extension is given by the middle term $\mu(B(1))$ of the
Auslander-Reiten sequence ending in $B(1)$. We refer to the end of section 
\ref{AR} for a closer look at $\mu(B(1)).$

\subsection{The use of the universal cover $\widetilde K(n)$ of $K(n)$}
We have used the universal cover $\widetilde K(n)$ of $K(n)$
in order to exhibit bristle submodules of some $n$-Kronecker modules $M$. 

{\bf (1)} Any leaf $x $ in the support of the covering module yields the bristle
module $P(x)$ for $\widetilde K(n)$, provided that $x$ is a source.
But there may be other bristle submodules of $M$ of the form $\pi(U)$, where $U$ is
a $\widetilde K(n)$-bristle. Let us look at the covering module $M$ for $I_3$ and,
on the left, its support quiver:
$$
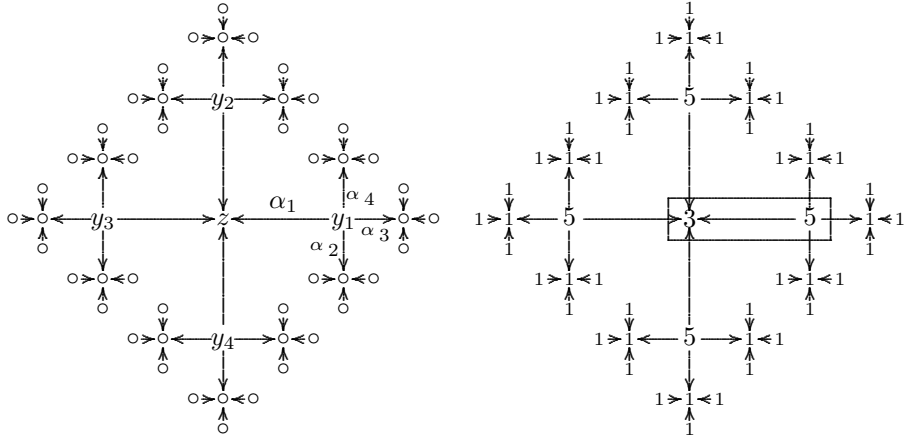
\hbox{\beginpicture
  \setcoordinatesystem units <.4cm,.4cm>
\put{\beginpicture

\put{$z$} at 0 0

\put{$y_1\strut$} at 4 0
\put{$\circ$} at 6 0
\put{$\circ$} at 7 0

\put{$\circ$} at 4 2
\multiput{$\circ$} at 4 3  3 2  5 2 /

\put{$\circ$} at 6 1

\arr{3.5 0}{0.3 0}
\arr{4.4 0}{5.7 0}
\arr{6.7 0}{6.3 0}

\arr{6 0.7}{6 0.3}

\arr{4 0.4}{4 1.7}
\arr{3.3 2}{3.7 2}
\arr{4 2.7}{4 2.3}
\arr{4.7 2}{4.3 2}

\put{$\circ$} at 4 -2
\multiput{$\circ$} at 4 -3  3 -2  5 -2 /

\put{$\circ$} at 6 -1
\arr{6 -.7}{6 -.3}

\arr{4 -.4}{4 -1.7}
\arr{3.3 -2}{3.7 -2}
\arr{4 -2.7}{4 -2.3}
\arr{4.7 -2}{4.3 -2}


\put{$y_3\strut$} at -4 0
\put{$\circ$} at -6 0
\put{$\circ$} at -7 0

\put{$\circ$} at -4 2
\multiput{$\circ$} at -4 3  -3 2  -5 2 /

\put{$\circ$} at -6 1

\arr{-3.5 0}{-0.3 0}
\arr{-4.4 0}{-5.7 0}
\arr{-6.7 0}{-6.3 0}

\arr{-6 0.7}{-6 0.3}

\arr{-4 0.4}{-4 1.7}
\arr{-3.3 2}{-3.7 2}
\arr{-4 2.7}{-4 2.3}
\arr{-4.7 2}{-4.3 2}

\put{$\circ$} at -4 -2
\multiput{$\circ$} at -4 -3  -3 -2  -5 -2 /

\put{$\circ$} at -6 -1
\arr{-6 -.7}{-6 -.3}

\arr{-4 -.4}{-4 -1.7}
\arr{-3.3 -2}{-3.7 -2}
\arr{-4 -2.7}{-4 -2.3}
\arr{-4.7 -2}{-4.3 -2}


\put{$y_2\strut$} at 0 4
\put{$\circ$} at 0 6
\put{$\circ$} at 0 7

\put{$\circ$} at 2 4
\multiput{$\circ$} at  3 4  2 3  2 5  /

\put{$\circ$} at 1 6

\arr{0 3.5}{0 0.3}
\arr{0 4.4}{0 5.7}
\arr{0 6.7}{0 6.3}

\arr{0.7 6}{0.3 6}

\arr{0.4 4}{1.7 4}
\arr{2 3.3}{2 3.7}
\arr{2.7 4}{2.3 4}
\arr{2 4.7}{2 4.3}

\put{$\circ$} at -2 4
\multiput{$\circ$} at -3 4   -2 3   -2 5 /

\put{$\circ$} at -1 6
\arr{-.7 6}{-.3 6}

\arr{-.4 4}{-1.7 4}
\arr{-2 3.3}{-2 3.7}
\arr{-2.7 4}{-2.3 4}
\arr{-2 4.7}{-2 4.3}


\put{$y_4\strut$} at 0 -4
\put{$\circ$} at 0 -6
\put{$\circ$} at 0 -7

\put{$\circ$} at 2 -4
\multiput{$\circ$} at  3 -4  2 -3  2 -5  /

\put{$\circ$} at 1 -6

\arr{0 -3.5}{0 -0.3}
\arr{0 -4.4}{0 -5.7}
\arr{0 -6.7}{0 -6.3}

\arr{0.7 -6}{0.3 -6}

\arr{0.4 -4}{1.7 -4}
\arr{2 -3.3}{2 -3.7}
\arr{2.7 -4}{2.3 -4}
\arr{2 -4.7}{2 -4.3}

\put{$\circ$} at -2 -4
\multiput{$\circ$} at -3 -4   -2 -3   -2 -5 /

\put{$\circ$} at -1 -6
\arr{-.7 -6}{-.3 -6}

\arr{-.4 -4}{-1.7 -4}
\arr{-2 -3.3}{-2 -3.7}
\arr{-2.7 -4}{-2.3 -4}
\arr{-2 -4.7}{-2 -4.3}

\put{$ \alpha_1$} at 2 0.5

\multiput{$\ssize \alpha_2$} at 3.3 -1  /
\multiput{$\ssize \alpha_3$} at 5 -.5   /
\multiput{$\ssize \alpha_4$} at 4.5 .7  /

\endpicture} at 0 0

\put{\beginpicture

\put{$3$} at 0 0

\put{$5\strut$} at 4 0
\put{$\ssize 1$} at 6 0
\put{$\ssize 1$} at 7 0

\put{$\ssize 1$} at 4 2
\multiput{$\ssize 1$} at 4 3  3 2  5 2 /

\put{$\ssize 1$} at 6 1

\arr{3.5 0}{0.3 0}
\arr{4.4 0}{5.7 0}
\arr{6.7 0}{6.3 0}

\arr{6 0.7}{6 0.3}

\arr{4 0.4}{4 1.7}
\arr{3.3 2}{3.7 2}
\arr{4 2.7}{4 2.3}
\arr{4.7 2}{4.3 2}

\put{$\ssize 1$} at 4 -2
\multiput{$\ssize 1$} at 4 -3  3 -2  5 -2 /

\put{$\ssize 1$} at 6 -1
\arr{6 -.7}{6 -.3}

\arr{4 -.4}{4 -1.7}
\arr{3.3 -2}{3.7 -2}
\arr{4 -2.7}{4 -2.3}
\arr{4.7 -2}{4.3 -2}


\put{$5\strut$} at -4 0
\put{$\ssize 1$} at -6 0
\put{$\ssize 1$} at -7 0

\put{$\ssize 1$} at -4 2
\multiput{$\ssize 1$} at -4 3  -3 2  -5 2 /

\put{$\ssize 1$} at -6 1

\arr{-3.5 0}{-0.3 0}
\arr{-4.4 0}{-5.7 0}
\arr{-6.7 0}{-6.3 0}

\arr{-6 0.7}{-6 0.3}

\arr{-4 0.4}{-4 1.7}
\arr{-3.3 2}{-3.7 2}
\arr{-4 2.7}{-4 2.3}
\arr{-4.7 2}{-4.3 2}

\put{$\ssize 1$} at -4 -2
\multiput{$\ssize 1$} at -4 -3  -3 -2  -5 -2 /

\put{$\ssize 1$} at -6 -1
\arr{-6 -.7}{-6 -.3}

\arr{-4 -.4}{-4 -1.7}
\arr{-3.3 -2}{-3.7 -2}
\arr{-4 -2.7}{-4 -2.3}
\arr{-4.7 -2}{-4.3 -2}


\put{$5\strut$} at 0 4
\put{$\ssize 1$} at 0 6
\put{$\ssize 1$} at 0 7

\put{$\ssize 1$} at 2 4
\multiput{$\ssize 1$} at  3 4  2 3  2 5  /

\put{$\ssize 1$} at 1 6

\arr{0 3.5}{0 0.3}
\arr{0 4.4}{0 5.7}
\arr{0 6.7}{0 6.3}

\arr{0.7 6}{0.3 6}

\arr{0.4 4}{1.7 4}
\arr{2 3.3}{2 3.7}
\arr{2.7 4}{2.3 4}
\arr{2 4.7}{2 4.3}

\put{$\ssize 1$} at -2 4
\multiput{$\ssize 1$} at -3 4   -2 3   -2 5 /

\put{$\ssize 1$} at -1 6
\arr{-.7 6}{-.3 6}

\arr{-.4 4}{-1.7 4}
\arr{-2 3.3}{-2 3.7}
\arr{-2.7 4}{-2.3 4}
\arr{-2 4.7}{-2 4.3}


\put{$5\strut$} at 0 -4
\put{$\ssize 1$} at 0 -6
\put{$\ssize 1$} at 0 -7

\put{$\ssize 1$} at 2 -4
\multiput{$\ssize 1$} at  3 -4  2 -3  2 -5  /

\put{$\ssize 1$} at 1 -6

\arr{0 -3.5}{0 -0.3}
\arr{0 -4.4}{0 -5.7}
\arr{0 -6.7}{0 -6.3}

\arr{0.7 -6}{0.3 -6}

\arr{0.4 -4}{1.7 -4}
\arr{2 -3.3}{2 -3.7}
\arr{2.7 -4}{2.3 -4}
\arr{2 -4.7}{2 -4.3}

\put{$\ssize 1$} at -2 -4
\multiput{$\ssize 1$} at -3 -4   -2 -3   -2 -5 /

\put{$\ssize 1$} at -1 -6
\arr{-.7 -6}{-.3 -6}

\arr{-.4 -4}{-1.7 -4}
\arr{-2 -3.3}{-2 -3.7}
\arr{-2.7 -4}{-2.3 -4}
\arr{-2 -4.7}{-2 -4.3}

\plot -.7 -.7  4.7 -.7  4.7 .7  -.7 .7  -.7 -.7 /

\endpicture} at 15.5 0
\endpicture}
$$
Consider in the subspace $M_{y_1}$ the intersection $A$ of the three maps with labels
$\alpha_1, \alpha_2,\alpha_3.$  This subspace $A$ is mapped under $\alpha_1$
injectively to $M_{z}$, thus it generates a direct sum of two bristles of type $B(1)$
(see also (B.4) in the appendix).

Since $\dim\Hom(B(1),I_3) = 11$, we see that we obtain sufficiently many bristles
of type $B(1)$ by looking at the 9 leaves of type $1$ and two of the bristles
generated by elements of $A$.

{\bf (2)} The universal cover $X$
of $I_2$ displays all the bristles of type $B(i)$ of $I_2$ (with $1\le i \le n$)
in a nice way:
as coming from bristles of $\widetilde I_2$ given by leaves in $Q$.
Namely, the number of leaves $x$ in $Q$
which are starting points of arrows with label $\alpha_i$ is $n-1$,
and, as we have mentioned already in \ref{book},  we have
$\dim\Hom(B(i),I_2) =  n-1.$

{\bf (3)} The universal cover $X$
of $I_2$ allows also to trace the bristles of type
$$
 B(1,2,\dots,n) = (k,k;1,1,\dots,1).
$$
For $1\le j \le n$, let $Y(j)$ be the sum of the submodules
$P(x)$ of $X$, where $x$ is a leaf and a neighbor of $y_j$. Then $\pi(Y(j))$ has
dimension vector $\ssize{\binom {n-1}1}$, and $X/\sum Y(j)$ is a direct sum of
$n-1$ copies of $S(z)$, and canonically isomorphic to $X_z$. On the one hand, we have
$\dim\Hom(B(1,2,\dots,n),I_2) = n-1$, on the other hand $\Hom(B(1,2,\dots,n),Y(j)) = 0$
for all $j$. This shows that $X$ has a submodule $U$ isomorphic to $B(1,2,\dots,n)^{n-1}$
such that the composition of the canonical maps $U \to X \to X/\sum Y(j)$
is surjective.


\section{Further considerations.}\label{related} 
		 				                      \medskip

\subsection{}\label{AR}
{\bf Proposition.} {\it Let $M$ be an indecomposable bristled module, not a bristle and not simple
projective. If $\tau M$ is a bristled module, then
the middle term $\mu(M)$ of the Auslander-Reiten sequence ending in $M$ is
a bristled module.

If $n\ge 2$, then the middle term $\mu(B)$ of the Auslander-Reiten sequence ending in a bristle $B$ is
not a bristled module.}

\begin{proof} Let $M$ be an indecomposable bristled module, not a bristle and not simple
projective. Then $M$ is not projective, thus we can consider
the middle term $\mu(M)$ of the Auslander-Reiten sequence ending in $M$. Since the canonical map
$\mu(M) \to M$ is right almost split,
all the maps $B \to M$ with $B$ a bristle can be lifted to $\mu(M)$.
Together with the maps $B \to \tau(M)$, where $B$ is a bristle, we obtain enough maps
from bristles to $\mu(M)$ in order to see that $\mu(M)$ is generated by bristles.

Now consider $\mu(B)$, where $B$ is a bristle and let $\epsilon\colon\mu(B) \to B$
be the canonical map. Assume that $\mu(B)$ is bristled. Then $\mu(B)$ is a factor module
of a direct sum of bristles, thus there is a map $f\colon B' \to \mu(B)$ with $B'$ a bristle
such that $\epsilon f \neq 0$. But any non-zero map between bristles is an isomorphism, thus
$\epsilon f$ is an isomorphism and therefore $\epsilon$ is a split epimorphism.
This is impossible, thus $\mu(B)$ cannot be bristled. 
\end{proof}

Let us have a closer look at $\mu(B)$, say for $B = B(1)$. It can be realized
as a submodule of $I_2$, as follows. As in section  \ref{I-two}, we consider the
universal cover $\widetilde K(n)$ of $K(n)$ with push-down functor $\pi$ and denote by
$X$ a representation of $\widetilde K(n)$ with $\pi(X) = I_2.$ 
Let $X''$ be the submodule of $X$ generated by the submodules $Y(j)$ with $2\le j \le n$ 
and the subspace $X_z$. Thus we have
$$
 X' \subset X'' \subset X,
$$
with $X'$ as defined in \ref{reg} (and $\pi(X') = \tau B(1)$).
In order to see that $\pi(X'') = \mu(B(1))$, we only have to observe 
that $X''$ is indecomposable and 
that $X''/X'$ is an indecomposable length 2 module whose support
is an arrow with label $\alpha_1$, thus $\pi(X''/X')$ is isomorphic to $B(1)$. 

For example in case $n = 4$, the dimension vector of $X''$ looks as follows:
$$
\hbox{\beginpicture
\setcoordinatesystem units <.8cm,.8cm>
\put{\beginpicture
\put{} at 0 3.5
\put{$3$} at 0 0
\multiput{$1$} at 2 0  0 2  -2 0  0 -2
    1 2  0 3  -1 2
    -2 1  -2 -1  -3 0
    1 -2  0 -3  -1 -2 /
\arr{0.3 0}{1.7 0}
\arr{-.3 0}{-1.7 0}
\arr{0 0.3}{0 1.7}
\arr{0 -.3}{0 -1.7}


\arr{-2.7 0}{-2.3 0}
\arr{-2 0.7}{-2 0.3}
\arr{-2 -.7}{-2 -.3}

\arr{0 2.7}{0 2.3}
\arr{0.7 2}{0.3 2}
\arr{-.7 2}{-.3 2}

\arr{0 -2.7}{0 -2.3}
\arr{0.7 -2}{0.3 -2}
\arr{-.7 -2}{-.3 -2}
\put{$ \alpha_1$} at 1 0.3
\put{$ \alpha_2$} at -.4 .95
\put{$ \alpha_3$} at -1 -.3
\put{$ \alpha_4$} at .4 -1

\multiput{$\ssize \alpha_1$} at -.5 2.25  -2.5 0.25  -.5 -1.75 /
\multiput{$\ssize \alpha_2$} at -.3 -2.6  -2.3  -.6 /
\multiput{$\ssize \alpha_3$} at  0.6 1.7  0.5 -2.3 /
\multiput{$\ssize \alpha_4$} at .35 2.55  -1.65 0.55 /
\put{$X''$} at -3 2.5
\endpicture} at 8 0
\endpicture}
$$
\vfill\eject

\subsection{} Let $n\ge 3.$ 
The Auslander-Reiten component containing a bristle $B$ has the following shape
$$
  
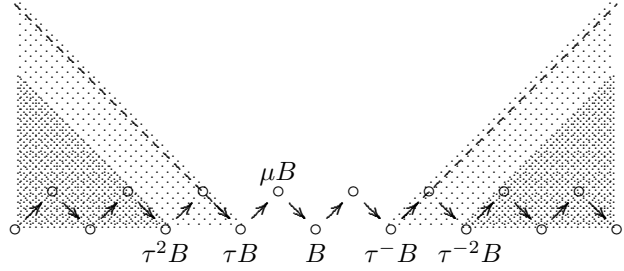
\hbox{\beginpicture
  \setcoordinatesystem units <.5cm,.5cm>
  \multiput{} at 0 0  0 6 /
  \multiput{$\circ$} at 0 0  2 0  4 0  6 0  8 0  10 0  12 0  14 0  16 0
            1 1  3 1  5 1  7 1  9 1  11 1  13 1  15 1 /
  \arr{0.3 0.3}{0.7 0.7}
  \arr{2.3 0.3}{2.7 0.7}
    \arr{4.3 0.3}{4.7 0.7}
  \arr{6.3 0.3}{6.7 0.7}
  \arr{8.3 0.3}{8.7 0.7}
  \arr{10.3 0.3}{10.7 0.7}
  \arr{12.3 0.3}{12.7 0.7}
  \arr{14.3 0.3}{14.7 0.7}
  \arr{1.3 0.7}{1.7 0.3}
  \arr{3.3 0.7}{3.7 0.3}
  \arr{5.3 0.7}{5.7 0.3}
  \arr{7.3 0.7}{7.7 0.3}
  \arr{9.3 0.7}{9.7 0.3}
  \arr{11.3 0.7}{11.7 0.3}
  \arr{13.3 0.7}{13.7 0.3}
  \arr{15.3 0.7}{15.7 0.3}
  \setdashes <1mm>
  \plot 0 6  5.7 0.3 /
  \plot 10.3 .3  16 6 /
  \put{$B$\strut} at 8 -.7
  \put{$\tau B$\strut} at 6 -.7
  \put{$\tau^2 B$\strut} at 4 -.7
  \put{$\mu B$\strut} at 7 1.5
  \put{$\tau^{-}B$\strut} at 10 -.7
 \put{$\tau^{-2}B$\strut} at 12.1 -.7
  \setshadegrid span <.7mm>
  \vshade 0 0 6.3  6.3 0 0 /
  \vshade 9.7 0 0  16 0 6.3 /
  \setshadegrid span <.4mm>
  \vshade 0 0 4.3  4.3 0 0 /
 \vshade 11.7 0 0  16 0 4.3 /
 \endpicture}
$$
 
\noindent
On the left hand side, we see the shaded cocone of all modules with
a path in the Auslander-Reiten quiver ending in $\tau B$,
all these modules are bristled modules; the darker shading on the left
marks the region of the modules which are in addition $\mathcal B$-saturated
(the assertions follow from \ref{cor1} and \ref{cor3} and the closure under
extensions).

On the right hand side, we see the shaded cone of all modules with
a path in the Auslander-Reiten quiver starting at $\tau^- B$.
By duality, all these modules are cogenerated by bristles.
Here, the darker shading indicates that the corresponding modules $M$ are in
addition $\mathcal B$-cosaturated.
	\medskip 

The remaining modules $M$ in the component (those which are neither 
predecessors of $\tau B$ nor successors of $\tau^-B$) all have bristle submodules
as well as bristle factor modules (since the northeast arrows correspond to monomorphisms,
the southeast arrows to epimorphisms); thus, with the exception of $B$ itself, are 
neither bristled nor cobristled. 
$$
  
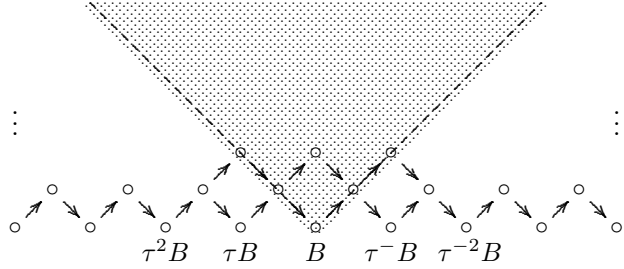
\hbox{\beginpicture
  \setcoordinatesystem units <.5cm,.5cm>
  \multiput{} at 0 0  0 6 /
  \multiput{$\circ$} at 0 0  2 0  4 0  6 0  8 0  10 0  12 0  14 0  16 0
            1 1  3 1  5 1  7 1  9 1  11 1  13 1  15 1 
                6 2  8 2  10 2  /
  \arr{0.3 0.3}{0.7 0.7}
  \arr{2.3 0.3}{2.7 0.7}
  \arr{4.3 0.3}{4.7 0.7}
  \arr{6.3 0.3}{6.7 0.7}
  \arr{8.3 0.3}{8.7 0.7}
  \arr{10.3 0.3}{10.7 0.7}
  \arr{12.3 0.3}{12.7 0.7}
  \arr{14.3 0.3}{14.7 0.7}
  \arr{1.3 0.7}{1.7 0.3}
  \arr{3.3 0.7}{3.7 0.3}
  \arr{5.3 0.7}{5.7 0.3}
  \arr{7.3 0.7}{7.7 0.3}
  \arr{9.3 0.7}{9.7 0.3}
  \arr{11.3 0.7}{11.7 0.3}
  \arr{13.3 0.7}{13.7 0.3}
  \arr{15.3 0.7}{15.7 0.3}

  \arr{6.3 1.7}{6.7 1.3}
  \arr{8.3 1.7}{8.7 1.3}
  \arr{10.3 1.7}{10.7 1.3}

  \arr{5.3 1.3}{5.7 1.7}
  \arr{7.3 1.3}{7.7 1.7}
  \arr{9.3 1.3}{9.7 1.7}

  \setdashes <1mm>

\plot 2 6  7.7 0.3 /
\plot 8.3 .3 14 6 /
  \put{$B$\strut} at 8 -.7
  \put{$\tau B$\strut} at 6 -.7
  \put{$\tau^2 B$\strut} at 4 -.7
  \put{$\tau^{-}B$\strut} at 10 -.7
 \put{$\tau^{-2}B$\strut} at 12.1 -.7
  \setshadegrid span <.5mm>
  \hshade -.2 7.8 8.2  <z,z,z,z> 6 1.8 14.2 /
\multiput{$\vdots$} at 0 3  16 3 /
 \endpicture}
$$

\subsection{}
{\it Let $n\ge 2.$
If $M$ is a module with a sectional path ending in a bristle $B$, then
$M$ is not $\mathcal B$-saturated.} (In particular, if
$\mu B$ is the middle term of the
Auslander-Reiten sequence ending in $B$, then $\mu B$ is not $\mathcal B$-saturated.)

\begin{proof} The epimorphism $M \to B$ yields a surjection
$\Ext^1(B,M) \to  \Ext^1(B, B).$
Calculating  $\langle \bdim B,\bdim B\rangle$, we see that
$\dim \Ext^1(B, B) = n-1$.
\end{proof}

\subsection{}
Actually, it is easy to see that for any bristle $B$, we have
$\dim \Ext^1(B,\mu B) = n-1.$ Namely, the Auslander-Reiten sequence ending in $B$ yields an exact sequence
$$
  \Hom(B,\mu B) \to \Hom(B,B) \to \Ext^1(\tau B,B) \to \Ext^1(\mu B,B) \to \Ext^1(B,B)
  	                          \to 0.
							        $$
The first map cannot be surjective, since otherwise the canonical projection $\mu B \to B$
would split. Since $\End(B) = k$, we see that the first map is the zero map, thus
the second map is injective. Again using that $\End(B) = k$, we
have $\dim \Ext^1(\tau B,B) = 1$. Thus, the second map is an isomorphism and therefore
the map $\Ext^1(\mu B,B) \to \Ext^1(B,B)$ is bijective.
As we have mentioned in the proof above, one may use the bilinear form
$\langle-,-\rangle$
in order to see that $\dim\Ext^1(B,B) = n-1$. Thus, also $\dim\Ext^1(\mu B,B) = n-1.$

The following picture shows for $n=4$ and $B = B(1)$ an indecomposable
representation $X'''$ of the universal cover $\widetilde K(4)$ with submodule
$X''$ such that $\pi(X'') = \mu (B(1)$) and
$\pi(X'''/X'')$ is isomorphic to $B(1)^3.$ Thus, the universal cover
exhibits nicely a basis of $\Ext^1(B(1),\mu(B(1))).$
$$
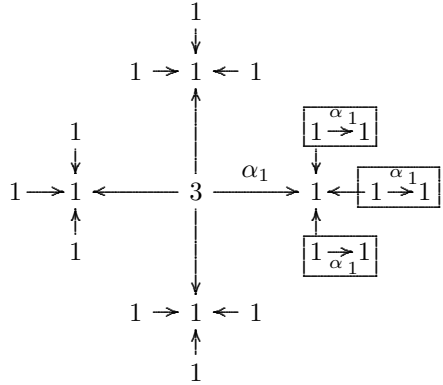
\hbox{\beginpicture
  \setcoordinatesystem units <.8cm,.8cm>
  \put{$3$} at 0 0
\multiput{$1$} at 2 0  0 2  -2 0  0 -2
    1 2  0 3  -1 2
    2 1  2 -1  3 0
    -2 1  -2 -1  -3 0
    1 -2  0 -3  -1 -2 /
\arr{0.3 0}{1.7 0}
\arr{-.3 0}{-1.7 0}
\arr{0 0.3}{0 1.7}
\arr{0 -.3}{0 -1.7}

\arr{2.8 0}{2.2 0}
\arr{2 0.7}{2 0.3}
\arr{2 -.7}{2 -.3}

\arr{-2.8 0}{-2.2 0}
\arr{-2 0.7}{-2 0.3}
\arr{-2 -.7}{-2 -.3}

\arr{0 2.7}{0 2.3}
\arr{0.7 2}{0.3 2}
\arr{-.7 2}{-.3 2}

\arr{0 -2.7}{0 -2.3}
\arr{0.7 -2}{0.3 -2}
\arr{-.7 -2}{-.3 -2}
\put{$ \alpha_1$} at 1 0.3


\arr{2.2 1}{2.6 1}
\arr{2.2 -1}{2.6 -1}
\arr{3.2 0}{3.6 0}
\multiput{$1$} at 2.8 1  2.8 -1  3.8 0 /
\multiput{$\ssize \alpha_1$} at 2.43 1.25  2.43 -1.25  3.43 0.25 /

\plot 1.8 0.75  3.0 0.75  3.0 1.4  1.8 1.4  1.8 0.75 /
\plot 2.7 -.25  4.05 -.25  4.05 .4  2.7 .4  2.7 -.25 /
\plot 1.8 -.75  3.0 -.75  3.0 -1.4  1.8 -1.4  1.8 -.75 /
\endpicture}
$$

We have just seen that $\mu(B(1))$ is not $\mathcal B$-saturated.
The following proposition shows that there are many modules, even many bristled modules,
which are not $\mathcal B$-saturated.

\subsection{}{\bf Proposition.}
{\it  An indecomposable $\mathcal B$-saturated module
is simple or faithful.}

\begin{proof}
Proof. For $n = 1$, all indecomposable modules are simple or faithful. Thus,
we may assume that $n\ge 2$. Let $M$ be
a non-faithful indecomposable $n$-Kronecker module.
Let as assume that $M$ is annihilated by
$\alpha_n$. We also assume that $M$ is different from $S(1)$. The following lemma
shows that $M$ is not $\mathcal B$-saturated.
\end{proof}

\subsection{Lemma.}{\it If $M$ is a $K(n)$-module which is annihilated by $\alpha_n$, then}
$$
 \dim\Ext^1(B(1),M) \ge \dim M_2.
 $$

\begin{proof}
We may consider both $B(1)$ and $M$ also as $(n-1)$-Kronecker modules.
We calculate $\langle \bdim B(1),\bdim M,\rangle_t$ for $t = n$ and for $t=n-1$.
We have
\begin{eqnarray*}
      \langle \bdim B(1),\bdim M\rangle_n     &=& \dim M_1-(n-1)\dim M_2, \\
            \langle \bdim B(1),\bdim M\rangle_{n-1} &=& \dim M_1-(n-2)\dim M_2,
\end{eqnarray*}
thus the difference is
$$
  D = \langle \bdim B(1),\bdim M\rangle_{n-1} -
          \langle \bdim B(1),\bdim M\rangle_{n} = \dim M_2.
$$
Since there is no difference whether we calculate $\Hom$ for $K(n)$ or for $K(n-1)$), 
the difference $D$ is also 
$$
  D =  -\dim\Ext^1_{K(n-1)}(B(1),M))+
                \dim\Ext^1_{K(n)}(B(1),M).
$$
Altogether, we see that
$$
    \dim\Ext^1_{K(n)}(B(1),M)) = \dim\Ext^1_{K(n-1)}(B(1),M)) + \dim M_2 \ge  \dim M_2.
$$
\end{proof}

\section*{Appendix A. The bristle variety $\beta(M)$ of a Kronecker module $M$.}

We assume here that $k$ is algebraically closed.

Let $M = (M_1,M_2,\alpha_1,\dots,\alpha_n)$
be a Kronecker module. First, we deal with the case where $M$ has no direct summand
of the form $S(1)$, thus the intersection of the kernels of the maps $\alpha_i$ is zero.
Let $\mathbb PM_1$ be the projective space corresponding to $M_1$, this is the set of all
1-dimensional subspaces of $M_1$, endowed with the Zariski topology. If $m$ is a non-zero
element of $M_1$, we write $\langle m\rangle$ for the $k$-subspace generated by $m$
(this is an element of $\mathbb PM_1$). Let $\beta(M)$ be
the subset of $\mathbb PM_1$ of all $\langle u\rangle$ such that $u$ generates a
bristle submodule in $M_1$ (or, equivalently, such that
$\alpha_1(u),\dots,\alpha_n(u)$ generate a one-dimensional subspace of $M_2$).
It is easy to see that $\beta(M)$ is closed in the Zariski topology, thus
{\it $\beta(M)$ is a projective variety,} we call it the {\it bristle variety}
of $M$. 

For a general $n$-Kronecker module $M$, the {\it bristle variety} of $M$ is by
definition $\beta(M/M'),$ where
$M'$ is the largest submodule of $M$ generated by $S(1)$ (so that $M/M'$ has no
submodule of the form $S(1)$).
	          \medskip

{\bf (A.1) Theorem.} {\it Any projective variety occurs as a bristle variety.}
            \medskip

The proof of Theorem (A.1) is given in \cite{[R1]}.
Note that this theorem is a variant of a sequence of results by Zimmermann-Huisgen, Hille,
Reineke and Van den Bergh which concern the realization of projective varieties say as
quiver Grassmannians. 
As a consequence of Theorem (A.1), the following is shown in \cite{[R2]}: {\it given any
connected wild acyclic quiver $Q$ with at least 3 vertices, then any projective variety occurs
as a quiver Grassmannian for a suitable representation of $Q$.}
   \medskip

A Kronecker module $M = (M_1,M_2,\alpha_1,\dots,\alpha_n)$ is nothing else than a
matrix pencil (as soon as we choose bases of the vector spaces $M_1,M_2$), and
generators of bristle submodules can be considered as (generalized) eigenvectors.
Indeed, if $n = 2$ and $\alpha_1$ is the identity map of a vector space
$V = M_1 = M_2$,
then a submodule of the $2$-Kronecker module $(V,V;1,\phi)$ is a bristle of type
$B_c = B(1:c)$ if and only if any generator of $U$ in $M_1$ is an
eigenvector of $\phi$ with eigenvalue $c$.
	                                       \medskip

\section*{Appendix B. Bristled modules for arbitrary artin algebras.}
     \medskip
     
Let $\Lambda$ be an artin algebra. Given a $\Lambda$-module $M$, let $|M|$ be its length.
Here, we collect some properties of the bristled $\Lambda$-modules and show in which way
the study of bristled $\Lambda$-modules can be reduced to the study of bristled Kronecker
modules.

	\medskip

{\bf (B.1) Proposition.} {\it If $M$ is an indecomposable bristled $\Lambda$-module $M$ 
then $|\topp M| \ge |\soc M|$, with equality only in case $M$ is simple or a bristle.} 
	\medskip 

\begin{proof} A bristled module has Loewy-length at most $2$. It follows that any
indecomposable bristled $\Lambda$-module $M$ is either simple or $\soc M = \rad M$.
Of course, if $M$ is simple, then $|\topp M| = 1 = |\soc M|.$ 

Thus, let us assume that $M$ is a bristled module with $\soc M = \rad M$. 
If $\topp M$ has length $t$, there is a
surjective map $f\colon \bigoplus_{i=1}^t B_i \to M$ with bristles $B_i$. 
The kernel of $f$ has to be contained in the socle of $\bigoplus B_i$,
let $s$ be its length. Then $|\soc M| = t-s \le t,$ this is the first assertion.
If $|\soc M| = t,$ then $s = 0$, thus $f$ is an isomorphism. If $M$ is indecomposable,
then we must have that $t=1$, thus $M$ is a bristle. 
\end{proof}

Recall that a semi-simple $\Lambda$-module
is said to be {\it homogeneous} provided that all its simple submodules are isomorphic.
   \medskip
   
{\bf (B.2) Proposition.} {\it The socle of an indecomposable bristled module $M$ is
   homogeneous.}

\begin{proof} Let $M$ be an indecomposable bristled module. Note that $M$ has Loewy length at most $2$, since $M$ is bristled. If $M$ has Loewy length $1$, then $M$ is
simple, therefore nothing has to be shown.
Thus, we assume that $M$ has Loewy length $2$. This means that
      $\rad M = \soc M$.

Let $S_1,\dots,S_n$ be the simple $\Lambda$-modules and write
$\soc M = \bigoplus_{i=1}^n H_i,$
where $H_i$ is the sum of all submodules of $M$ isomorphic to $S_i$.

By assumption, $M$ is generated by bristles, thus $M$ is the sum of images
of some maps $B\to M$, where $B$ is a bristle. Since $\rad M = \soc M$, it follows that
$M$ is generated by the images of some maps $B\to M$ which are in addition injective, thus
$M$ is the sum of submodules $U_j$ which are isomorphic to bristles, say
$M = \sum_{j\in J} U_j$. We may assume that $J$ is minimal, or, equivalently, that
$M/\soc M = \bigoplus (U_j+\soc M)/\soc M.$

Let $M_i$ be the sum of the submodules $U_j$ such that $\soc U_j$ is isomorphic to $S_i$.
Now $\soc U_j$ is a submodule of $\soc M$ and isomorphic to $S_i$, thus
$\soc U_j \subseteq H_i$. It follows easily that $M = \bigoplus_{i=1}^n M_i$.
Since we assume
that $M$ is indecomposable, we have $M = M_i$ for some $i$.
Of course, the socle of $M_i$ is just $H_i$, thus homogeneous.
\end{proof}

{\bf (B.3)} We say that a module $M$ is the {\it socle amalgamation of the modules
$M_1,\dots, M_t$}
provided the modules $M_1,\dots, M_t$ are submodules of $M$ containing the
socle $\soc M$ of $M$ and $M/\soc M$ is the direct sum $\bigoplus_i M_i/\soc M$.

It follows from (B.2) that {\it any indecomposable bristled module $M$
is the socle amalgamation of bristled modules $M_1,\dots,M_t$, such that
$\soc M$ as well all the modules $M_i/\soc M$ are homogeneous.}

Namely,
$\overline M = M/\soc M$ is
semi-simple, thus we may write $\overline M$
as a direct sum $\bigoplus_i \overline M_i$ of pairwise $\Hom$-orthogonal
homogeneous modules 
$\overline M_i$. Let $M_i$ be the preimage of $\overline M_i$ under
the canonical projection $M \to \overline M$.
    \medskip

     Thus, the first aim for understanding bristle modules seems to be to investigate
     bristled modules $M$ such that both $\soc M$ and $M/\soc M$
     are homogeneous.

    In case $\Lambda$ is a finite-dimensional $k$-algebra, where $k$
is an algebraically closed field,
it is sufficient to look just at the Kronecker algebras. Namely,
we may assume that $M$ is faithful and that $\Lambda$ is basic.
Under these assumptions, $\Lambda$ is a basic radical-square-zero algebra with
at most two simple modules. Since $k$ is algebraically closed, only two
cases are possible. If $\Lambda$ has two
simple modules, then $\Lambda$ is already a Kronecker algebra.
If $\Lambda$ has only one simple module, then $\Lambda$ is
of the form $k[x_1,\dots,x_n]/(x_1,\dots,x_n)^2,$ and therefore
stably equivalent to the Kronecker algebra $kK(n)$. Under such a
stable equivalence,
the indecomposable bristled $\Lambda$-modules correspond bijectively to the
indecomposable bristled $kK(n)$-modules.
            \bigskip

{\bf (B.4)} Let us assume now that $\Lambda$ is the path algebra of a
finite acyclic quiver $Q$ without multiple arrows. In this case, the support
of a bristle $B$ is an arrow. Attaching to any arrow $\alpha$ the thin indecomposable
module $B(\alpha)$ with support $\alpha$, we obtain a bijection between the arrows and
the (isomorphism classes of the) bristles. There is an easy numerical criterion to
detect (at least some) bristle submodules of a $kQ$-module $M$:
If $\alpha\colon  y \to z$ is an arrow and
$z_1 = z, z_2,\dots,z_t$ are the terminal vertices of the arrows starting in $y$,
then $\dim\Hom(B(\alpha),M) \ge \dim M_y-\sum_{i\ge 2} \dim M_{z_i}.$
(For the proof, one just has to look at $\langle \bdim B(\alpha),\bdim M\rangle.$)
Of course, if $B$ has no submodule isomorphic to $S(y)$, then
$\dim\Hom(B(\alpha),M) = b$ implies that $M$ has a submodule isomorphic to $B(\alpha)^b.$

In order to see that in general not all bristles are detected in this way,
consider the following $\widetilde{\mathbb D}_5$-quiver shown below on the left.
There exists an indecomposable representation $M$ with dimension vector presented
below on the right, which has a submodule $B(\alpha).$
whereas $\dim M_y-\sum_{i\ge 2} \dim M_{z_i} = 0.$
$$
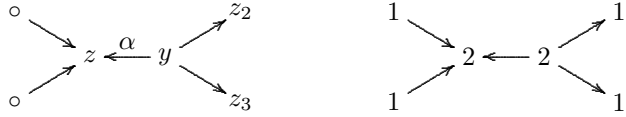
\hbox{\beginpicture
  \setcoordinatesystem units <1cm,.6cm>
  \put{\beginpicture
  \multiput{$\circ$} at 0 1  0 -1  /
  \put{$z$} at 1 0
  \put{$y$} at 2 0
  \put{$z_2$} at 3 1
  \put{$z_3$} at 3 -1
  \put{$\alpha$} at 1.5 0.3
  \arr{0.2 0.8}{0.8 0.2}
  \arr{0.2 -.8}{0.8 -.2}
  \arr{1.8 0}{1.2 0}
  \arr{2.2 0.2}{2.8 0.8}
  \arr{2.2 -.2}{2.8 -.8}
  \endpicture} at 0 0
  \put{\beginpicture
  \multiput{$1$} at 0 1  0 -1  /
  \put{$2$} at 1 0
  \put{$2$} at 2 0
  \put{$1$} at 3 1
  \put{$1$} at 3 -1
  \arr{0.2 0.8}{0.8 0.2}
  \arr{0.2 -.8}{0.8 -.2}
  \arr{1.8 0}{1.2 0}
  \arr{2.2 0.2}{2.8 0.8}
  \arr{2.2 -.2}{2.8 -.8}
  \endpicture} at 5 0
  \endpicture}
  $$
  \medskip

  {\bf (B.5)} Using the universal covering $\widetilde K(n)$ of $K(n)$, we have
  shown in the paper that all the preinjective $K(n)$-modules are bristled. But we should note
  that the non-injective preinjective $\widetilde K(n)$-modules themselves are not bristled.
       \medskip

{\bf Proposition.} {\it Let $Y$ be an indecomposable preinjective $\widetilde K(n)$-module.
If $Y$ is injective, then $Y$ is bristled. If $Y$ is not injective, then $Y$ is not bristled.}

\begin{proof} It is obvious that the indecomposable injective $\widetilde K(n)$-modules are
bristled. Thus, we only have to consider the remaining indecomposable preinjective
$\widetilde K(n)$-modules.

First, we consider the $\widetilde K(n)$-module $X = X(z) = \widetilde\tau S(z)$,
where $z$ is a source of $\widetilde K(n)$
(and $\widetilde\tau$ the Auslander-Reiten translation for
$\widetilde K(n)$); note that $X$ has been discussed already in \ref{I-two}. As before,
we denote by $\alpha_i\colon z \to y_i$ the arrows starting at $z$. The thin
indecomposable representation of $\widetilde K(n)$ with support $\alpha_i\colon z \to y_i$
will be denoted by $B_i$. Note that the modules $B_i$ are the only bristles with top
$S(z)$.

We claim that $\Hom(B_i,X) = 0.$ Denote by $Q'$ the full subquiver
of $\widetilde K(n)$ with vertices $z,y_1,\dots,y_n$. The restriction $X'$ of $X$ to $Q'$ is
indecomposable (and not injective). Since $B_i$ and all its factor modules are injective,
when considered as representations of $Q'$, we see that $\Hom(B_i,X') = 0,$ thus also
$\Hom(B_i,X) = 0.$

If $\mathcal B(X)$ is the maximal bristled submodule of $X$, then the top of $\mathcal B(X)$
has no composition factor of the form $S(z)$. As a consequence, $\mathcal B(X)$ is a proper
submodule of $X$. Thus $X$ is not bristled.

If $Y$ is an indecomposable preinjective representation of $\widetilde K(n)$ which is
not injective, then $Y$ has a factor module of the form $X(z)$ for some $z$. Since
$X(z)$ is not bristled, also $Y$ cannot be bristled.
\end{proof}
	  
\section*{Appendix C. A non-finiteness assertion for tame algebras.}
     \medskip
     
{\bf (C.1) Proposition.} {\it Let $\Lambda$ be a tame hereditary artin algebra.
     There is no regular module which generates infinitely many indecomposable
     modules. In particular, there is no regular or preinjective module which generates
     all preinjective modules.}

\begin{proof} Obviously, a preinjective module generates only finitely many indecomposable modules (and these are preinjective again).

If $R$ is simple regular and $t\in \mathbb N,$ then
the only regular modules generated by $R[t]$ are modules in the tube of $R$ of
regular length at most $t$, thus only finitely many isomorphism classes.
This shows that a regular module generates only finitely many indecomposable regular
modules.

Finally, let $R$ be indecomposable regular and $M$ indecomposable
preinjective. Then $\dim\Hom(R,M) \le rs$ where $r$ is the regular length of $R$
and $s$ is the maximal length of a simple regular module in the tube of $R$
(use $\tau$ in order to reduce to the
case that $M$ is indecomposable injective; under this shift, $r$ and $s$ do not change).
It follows that $M$ is a factor module of $R^{rs}$, thus of bounded length.
\end{proof}

{\bf (C.2)} 
In the case of the $2$-Kronecker quiver, one may strengthen the assertion as follows:
              \medskip

{\bf Proposition.} {\it Let $\Lambda = kK(2)$ and $M$ in $\mo\Lambda$. Then the
     		   	              following assertions are equivalent:}

(i) {\it $M$ generates infinitely many indecomposable modules.}

(ii) {\it $M$ generates all indecomposable modules with the exception of at most
   finitely many preprojective modules.}

(iii) {\it $M$ has an indecomposable preprojective direct summand.}
      \smallskip

      {\it If the preprojective module
      $P_t$ is a direct summand of $M$ and $t$ is minimal, then the only indecomposable
      modules which are not generated by $M$ are the modules $P_0,\dots, P_{t-1}.$}
      	                              \medskip

The proof is left to the reader.
    	     	              \hfill$\square$

 \bigskip
 
				        \bigskip

	\bigskip
	
\end{document}